\input amstex
\documentstyle{amsppt}
\NoBlackBoxes
\nopagenumbers
\magnification=\magstep1
\advance\hsize1cm\hoffset0cm
\advance\vsize0cm\voffset=-1.5cm
\define\R{{\Bbb R}}  \define\Z{{\Bbb Z}}
\define\C{{\Bbb C}} \define\Ha{{\Bbb H}}
\def\pr{\mathop{\fam0 pr}}
\def\rk{\mathop{\fam0 rk}}
\def\Sq{\mathop{\fam0 Sq}}

\def\Tors{\mathop{\fam0 Tors}}
\def\diag{\mathop{\fam0 diag}}

\def\Emb{\mathop{\fam0 E}}

\def\id{\mathop{\fam0 id}}

\def\lk{\mathop{\fam0 lk}}

\def\Int{\mathop{\fam0 Int}}

\def\Cl{\mathop{\fam0 Cl}}
\def\Emb{\mathop{\fam0 E}}

\def\im{\mathop{\fam0 im}}

\def\ct{1.1\ }
\def\Ct{1.1}

\def\Cor{1.2}

\def\adm{1.3\ }
\def\Adm{1.3}

\def\Coro{1.4}

\def\is{2.1\ }
\def\Is{2.1}

\def\ra{2.2}
\def\Rb{2.3\ }
\def\rb{2.3}
\def\nb{2.4\ }
\def\Nb{2.4}
\def\agr{2.5\ }
\def\Agr{2.5}
\def\co{2.6\ }
\def\Co{2.6}
\def\bo{2.7\ }
\def\Bo{2.7}
\def\hs{2.8\ }
\def\Hs{2.8}

\def\Ind{2.9}
\def\ad{2.10\ }
\def\Ad{2.10}
\def\di{2.11\ }
\def\Di{2.11}
\def\re{2.12\ }
\def\Re{2.12}
\def\eu{3.1\ }
\def\Eu{3.1}
\def\ide{3.2\ }
\def\Ide{3.2}
\def\ir{3.3\ }
\def\Ir{3.3}

\def\Kr{4.1}
\def\ho{4.2\ }
\def\Ho{4.2}
\def\bd{4.3\ }
\def\Bd{4.3}
\def\se{4.4\ }
\def\Se{4.4}
\def\tw{4.6\ }
\def\Tw{4.6}
\def\stw{4.5\ }
\def\Stw{4.5}

\def\Dif{4.7}

\def\Cjc{4.8}
\def\sh{4.9\ }
\def\Sh{4.9}
\def\li{5.5\ }
\def\Li{5.5}

\def\Cc{5.4}
\def\un{5.8\ }
\def\Un{5.8}


\topmatter
\title A classification of smooth embeddings of
4-manifolds in 7-space, II \endtitle
\author Diarmuid Crowley and Arkadiy Skopenkov \endauthor
\address
Hausdorff Research Institute for Mathematics, Universit\"at Bonn, Poppelsdorfer Allee 82, D-53115 Bonn, Germany. \vskip 0cm e-mail: diarmuidc23\@gmail.com
\endaddress
\address
Department of Differential Geometry, Faculty of Mechanics and
Mathematics, Moscow State University, 119992, Moscow, Russia, and
Independent University of Moscow, B. Vlasyevskiy, 11, 119002, Moscow, Russia.
\vskip 0cm e-mail: skopenko\@mccme.ru \endaddress
\abstract
Let $N$ be a closed connected smooth 4-manifold with $H_1(N;\Z)=0$.
Our main result is the following classification of the set $\Emb^7(N)$ of
smooth embeddings $N\to\R^7$ up to smooth isotopy.
Haefliger proved that $\Emb^7(S^4)$ together with the connected sum
operation is a group isomorphic to $\Z_{12}$.
This group acts on $\Emb^7(N)$ by embedded connected sum.
Bo\'echat and Haefliger constructed an invariant $\varkappa :\Emb^7(N)\to H_2(N;\Z)$
which is injective on the orbit space of this action;
they also described $\im\varkappa $.
We {\it determine the orbits} of the action:
{\it for  $u\in\im\varkappa $ the number of elements in $\varkappa ^{-1}(u)$ is
$GCD(u/2,12)$ if $u$ is divisible by $2$, or is $GCD(u,3)$ if $u$ is not
divisible by $2$.}  The proof is based on Kreck's modified formulation of surgery.
\endabstract
\endtopmatter

\document
\head 1. Introduction and main results  \endhead

The main result of this paper is {\it a complete readily calculable
classification of smooth embeddings into $\R^7$ of closed, smooth
4-manifolds $N$ such that $H_1(N)=0$}. Cf. [Sk10, footnote 1].
We work in the smooth category.
For a manifold $N$ let $\Emb^m(N)$ denote the set of smooth embeddings
$N\to\R^m$ up to smooth isotopy.
We omit $\Z$-coefficients from the notation of (co)ho\-mo\-lo\-gy groups and
denote Poincar\'e duality by $PD$.

We define the Bo\'{e}chat--Haefliger invariant and the Kreck invariant used in the following theorem in \S1 and \S2, respectively.

\smallskip
{\bf Classification Theorem \Ct.}
{\it Let $N$ be a closed connected 4-manifold such that $H_1(N)=0$.
The Bo\'echat-Haefliger invariant
$$\varkappa :\Emb\phantom{}^7(N)\to H_2(N)$$
has image
$$\im\varkappa =\{u\in H_2(N)\ |\ u\equiv PDw_2(N)\mod2,\ u\cap u=\sigma(N)\}.$$
For each $u\in\im\varkappa $ the Kreck invariant
$$\eta_u:\varkappa ^{-1}(u)\to\Z_{GCD(u,24)}$$
is injective and has image the subset of even elements.
\footnote{Here $GCD(u,24)$ is the maximal integer $k$ such that both
$u\in H_2(N)$ and 24 are divisible by $k$.
Thus $\eta_u$ is surjective if $u$ is not divisible by 2.
Note that $u\in\im\varkappa $ is divisible by 2 (for some $u$ or, equivalently, for
each $u$) if and only if $N$ is spin.}
}

\smallskip
{\bf Corollary \Cor.}
\footnote{For an explicit construction of the embeddings see \S3 and
Corollary \Coro(c) below.}
(a) {\it There are exactly twelve isotopy classes of
embeddings $N\to\R^7$ if $N=S^4$ [Ha66] or an integral homology 4-sphere.}


(b) {\it For each integer $u$ there are exactly $GCD(u,12)$ isotopy classes of
embeddings $f:S^2\times S^2\to\R^7$ with $\varkappa (f)=(2u,0)$, and
the same holds for those with $\varkappa (f)=(0,2u)$.
Other values of $\Z^2$ are not in the image of $\varkappa $.
(We take the standard basis in $H_2(S^2\times S^2)$.)}

\smallskip
The description of $\im\varkappa $ in the Classification Theorem \ct was already
known [BH70], cf. [Fu94].
So our achievement is {\it to describe the preimages of $\varkappa $} (thus only
this part of the proof is presented in this paper).
More precisely, in this description our achievement is the transition
from the case $N=S^4$ (which was known)
to closed connected 4-manifolds $N$ with $H_1(N)=0$.
\footnote{A simpler proof of a particular case of the Classification Theorem
\ct is given in [Sk10].}
Let us explain what is involved in this transition.

From now on unless otherwise stated, we assume the following:

$\bullet$ {\it $N$ is a closed connected orientable 4-manifold and $f:N\to\R^7$ is an
embedding}.

It was known that $\Emb^7(S^4)$ with the embedded connected sum operation
is a group isomorphic to $\Z_{12}$ [Ha66].
The group $\Emb^7(S^4)$ acts on the set
$\Emb^7(N)$ by connected summation of embeddings $g:S^4\to\R^7$ and
$f:N\to\R^7$ whose images are contained in disjoint cubes.
It was known that {\it for $H_1(N)=0$ the orbit space of this action
$\Emb^7(S^4)\to\Emb^7(N)$ maps bijectively under $\varkappa $ (defined in a different
way) to $\im\varkappa $. }
This follows by [BH70, Theorems 1.6 and 2.1] and
smoothing theory [BH70, p. 156].

\smallskip
{\bf Addendum \Adm.}
{\it Let $N$ be a closed connected 4-manifold such that $H_1(N)=0$.
For each pair of embeddings $f:N\to\R^7$ and $g:S^4\to\R^7$
$$\varkappa (f\#g)=\varkappa (f)\quad\text{and}\quad
\eta_{\varkappa (f)}(f\#g)\equiv\eta_{\varkappa (f)}(f)+\eta_0(g)\mod GCD(\varkappa (f),24).$$}
\quad
Here the first equality follows by the definition of the
Bo\'echat-Haefliger invariant, and the second equality is proved in \S3.

\smallskip
{\it Definition of the Bo\'echat-Haefliger invariant.}
Denote by $C_f$ the closure of the complement in $S^7\supset\R^7$ to a tubular
neighborhood of $f(N)$.  Fix an orientation on $N$ and an orientation on $\R^7$.
A {\it homology Seifert surface} $A_f$ for $f$ is the generator of
$H_5(C_f,\partial)\cong\Z$ chosen by the fixed orientations of $N$ and $\R^7$.
\footnote{More precisely, $A_f$ is the image of the fundamental class $[N]$ under
the composition $H_4(N)\to H^2(C_f)\to H_5(C_f,\partial)$ of the Alexander and
Poincar\'e-Lefschetz duality isomorphisms; this composition is the inverse to
the composition $H_5(C_f,\partial)\to H_4(\partial C_f)\to H_4(N)$ of the
boundary map and the normal bundle map, cf. [Sk08', the Alexander Duality
Lemma]; the latter assertion justifies the name `homology Seifert surface'.}

Define $\varkappa (f)$ to be the image of $A_f^2=A_f\cap A_f$ under the composition
$H_3(C_f,\partial)\to H^4(C_f)\to H_2(N)$ of the Poincar\'e-Lefschetz and
Alexander duality isomorphisms.

\smallskip
This new definition is equivalent to the original one [BH70] by
[Sk10, Remark 2.4, footnote 14 and the first equality of Section Lemma 2.5],
cf. Section Lemma 3.1 below.

The Classification Theorem \ct and Addendum \adm imply the following
{\it examples of the triviality and the effectiveness of the above action
$\Emb^7(S^4)\to \Emb^7(N)$}.


\smallskip
{\bf Corollary \Coro.}
(a) {\it Let $N$ be a closed connected 4-manifold such that $H_1(N)=0$ and
the signature $\sigma(N)$ of $N$ is not divisible by the square of an integer
$s\ge2$ (in particular, $N=\C P^2$).
Then for each embedding $f:N\to\R^7$ and $g:S^4\to\R^7$ the embedding $f\#g$
is isotopic to $f$ {\rm [Sk10, the Triviality Theorem 1.1.a].}
\footnote{In other words, under the assumption of Corollary \Coro(a) the map $\varkappa $ is
injective.}}

(b) {\it If $N$ is a closed connected 4-manifold such that $H_1(N)=0$ and
$f(N)\subset\R^6$ for an embedding $f:N\to\R^7$,
then for each embedding $g:S^4\to\R^7$ the embedding $f\# g$ is not isotopic to
$f$. Cf. {\rm [Sk10, the Effectiveness Theorem 1.2].}}

(c) {\it Take an integer $u$ and an embedding $f_u:S^2\times S^2\to\R^7$
constructed just below.
If $u=6k\pm1$, then for each embedding $g:S^4\to\R^7$ the embedding $f_u\#g$ is isotopic to
$f_u$.
\footnote{For a general integer $u$ the number of isotopy classes of embeddings
$f_u\#g$ is $GCD(u,12)$.}
}

\smallskip
{\it Sketch of a proof.}
Part (a) follows by Addendum \adm and the Classification Theorem \Ct.

Part (b) follows by the Classification Theorem \ct because $\varkappa (f)=0$ when
$f(N)\subset\R^6$, cf. [Sk08', Compression Theorem].

Part (c) follows by the Classification Theorem \ct because
$\varkappa (f_u)=2W(f_u)=2u$ analogously to [Sk08', Bo\'echat-Haefliger Invariant
Theorem], where $W(f_u)$ is defined analogously to [Sk08', definition of
the Whitney invariant].
\qed

\smallskip
{\it The first construction of $f_u$.}
Let $\overline f_u:S^2\to V_{5,3}$ be a map representing
$u$ times the generator of $\pi_2(V_{5,3})\cong\Z$.
This map $f_u$ can be seen as a map from $S^2$ to the space of linear
orthogonal embeddings $\R^3\to\R^5$.
By the exponential law this gives a map
$\widehat f_u=\pr_1\times\overline f_u:S^2\times\R^3\to S^2\times\R^5$, where
$\pr_1$ is the projection onto the first factor.
Let $f_u$ be the composition
$S^2\times\partial D^3\to S^2\times\partial D^5\to\R^7$ of
the restriction of $\widehat f_u$ and the standard inclusion.

\smallskip
{\it The second construction of $f_u$.}
Take the standard embeddings $2D^5\times S^2\subset\R^7$ (where 2 is
multiplication by 2) and $\partial D^3\subset\partial D^5$.
Take $u$ copies $(1+\frac1n)\partial D^5\times x$ ($n=1,\dots,u$) of the oriented 4-sphere
outside $D^5\times S^2$ `parallel' to $\partial D^5\times x$.
Join these spheres by tubes so that the homotopy class of the resulting
embedding $S^4\to S^7-(D^5\times S^2)\simeq S^7-S^2\simeq S^4$ will be
$u\in\pi_4(S^4)\cong\Z$.
Let $f$ be the connected sum of this embedding with the standard embedding
$\partial D^3\times S^2\subset\R^7$.


\smallskip
It follows from the Classification Theorem \ct that if $f_k:N_k\to\R^7$ are
embeddings of closed connected 4-manifolds such
that  $H_1(N_k)=0$ and $a_k:=\varkappa _{N_k}(f_k)$, then
$$\# \varkappa _{N_1\# N_2}^{-1}(a_1\oplus a_2)=
\cases GCD(a_1,a_2,3)& \text{if either $a_1$ or $a_2$ is not divisible by 2,}\\
GCD(a_1/2,a_2/2,12) &\text{if both $a_1$ and $a_2$ are divisible by 2.}
\endcases$$
\quad
We plan to prove a generalization of the Classification Theorem 1.1
to non-simply connected 4-manifolds in [CS].

\bigskip
{\bf The general Knotting Problem.}

This subsection gives some background about the Knotting Problem: it is not used in the proof of the Classification Theorem \Ct.  The classical Knotting Problem runs as follows: {\it given an $n$-manifold $N$ and a number $m$, describe $\Emb^m(N)$, the set of isotopy classes of embeddings $N\to\R^m$}.
\footnote{The classification of embeddings into $S^m$ is the same because
{\it if the compositions with the inclusion $i:\R^m\to S^m$ of two
embeddings $f_0,f_1:N\to\R^m$ of a compact $n$-manifold $N$ are isotopic, then
$f_0$ and $f_1$ are isotopic} (in spite of the existence of orientation-preserving
diffeomorphisms $S^m\to S^m$ not isotopic to the identity).
Indeed, since $f_0$ and $f_1$ are isotopic, by general position
$i\circ f_0$ and $i\circ f_1$ are non-ambiently isotopic.
Since every non-ambient isotopy extends to an ambient one [Hi76, Theorem 1.3],
$i\circ f_0$ and $i\circ f_1$ are isotopic.}
For recent surveys see [RS99, Sk08, MA2]; whenever possible we refer to these
surveys not to original papers.

The Knotting Problem is more accessible for $2m\ge3n+4$ [RS99, Sk08, \S2, \S3, MA2].
It is much harder for
$$2m<3n+4:$$
if $N$ is a closed manifold that is not a disjoint union of homology spheres, then
until recently no complete readily calculable descriptions of isotopy classes
was known, in spite of the existence of interesting approaches of Browder-Wall
and Goodwillie-Weiss [Wa70, GW99, CRS04].
\footnote{The approach of [GW99] gives a modern abstract proof of certain
earlier known results. We are grateful to M. Weiss for indicating that this
approach also gives explicit results on higher homotopy groups of the space of
embeddings $S^1\to\R^n$.}
For recent results see [Sk06, Sk08']; for {\it rational} and {\it piecewise
linear} classification see [CRS07, CRS08] and
[Sk06, Sk07, Sk08, \S2, \S3 and \S5], respectively.

In particular, a complete, readily calculable classification of embeddings of a
closed connected 4-manifold $N$ into $\R^m$ was only known only for $m\ge8$
(Wu, Haefliger, Hirsch and Bausum) or for $N=S^4$ and $m=7$ (Haefliger):
$$\#\Emb\phantom{}^m(N)=1\quad\text{for}\quad m\ge9,$$
$$\Emb\phantom{}^8(N)=\cases H_1(N;\Z_2) & N\text{ orientable,}\\
\Z\oplus\Z_2^{s-1} & N\text{ non-orientable and }H_1(N;\Z_2)\cong\Z_2^s,
\endcases$$
$$\Emb\phantom{}^7(S^4)\cong\Z_{12}.$$
Here the equality sign between sets denotes the existence of a bijection; the
isomorphism is a group isomorphism for certain geometrically defined group
structures.
See references and more information in [MA1].

The `connected sum' group structure on $\Emb^m(S^n)$ was defined in [Ha66].
By [Ha61, Ha66, Corollary 6.6, Sk08, \S3],
$\Emb\phantom{}^m(S^n)=0\quad\text{for}\quad 2m\ge3n+4.$
However, $\Emb\phantom{}^m(S^n)\ne0$ for many $m,n$ such that $2m<3n+4$,
\footnote{This differs from the Zeeman-Stal\-lings Unknotting Theorem:
{\it for $m\ge n+3$ any PL or TOP embedding $S^n\to S^m$ is PL or TOP
isotopic to the standard embedding}.}
e.g. $\Emb^7(S^4)\cong\Z_{12}$.

In this paragraph assume that $N$ is a closed $n$-manifold and $m\ge n+3$.
The group $\Emb^m(S^n)$ acts on the set
$\Emb^m(N)$ by connected summation of embeddings $g:S^n\to\R^m$ and
$f:N\to\R^m$ whose images are contained in disjoint cubes.
\footnote{Since $m\ge n+3$, the connected sum is well-defined, i.e. does not depend on
the choice of an arc between $gS^n$ and $fN$.
If $N$ is not connected, we assume that a component of $N$ is chosen and
we consider embedded connected summation with this chosen component. }
Various authors have studied the analogous connected sum action of the group of
homotopy $n$-spheres on the set of smooth $n$-manifolds homeomorphic to
given manifold; see for example [Le70].  For embeddings, the quotient of $\Emb^m(N)$ modulo the above action of $\Emb^m(S^n)$ is known in some cases.
\footnote{In those cases when this quotient coincides with $\Emb^m_{PL}(N)$
and when the latter set was known [Hu69, \S12, Vr77, Sk97, Sk02, Sk07, Sk06].}
Thus in these cases the knotting problem is reduced to {\it the determination
of the orbits of this action}.
This problem is just as difficult as the Knotting Problem: until recently {\it no results were known} on this action for $m\ge n+3$, $\Emb^m(S^n)\ne0$ and $N$ not a disjoint union of spheres.
For recent results see [Sk08', Sk06]; for a {\it rational} description see
[CRS07, CRS08]; for $m=n+2$ see [Vi73].

\bigskip
{\bf Acknowledgements.}
These results are based on ideas of and discussion with Matthias Kreck.
They were announced at the International Pontryagin Conference (Moscow, 2008).

\head 2. An overview of the proof \endhead

This section consists of four subsections.
The first discusses the general strategy we use.
The second states the preliminary results needed to apply this strategy to
calculate $\Emb^7(N)$.
The third defines the key invariant, the Kreck invariant.
The final subsection gives the proof of the Classification Theorem \ct
modulo some results proved in \S\S3-4.

\bigskip
{\bf A general strategy for the embedding problem.}
\nopagebreak

The proof of the Classification Theorem \ct is based on the ideas we
explain below which are useful in a wider range of dimensions [Sk08']
and for solving problems other than finding the action of $\Emb^m(S^n)$ on
$\Emb^m (N)$ [FKV87, FKV88].


In this subsection $N$ is a closed connected $n$-manifold and $f:N\to\R^m$ is an embedding.
Let $\nu_f$ be the normal vector bundle of $f(N)$ and let $C_f$ be the closure
of the complement in $S^m\supset\R^m$ to a tubular neighborhood of $f(N)$.
We identify the boundary of $C_f$, $\partial C_f$, with the total space
of the sphere bundle of $\nu_f$.
In this paper an oriented (or spin) bundle isomorphism is always the restriction of
an oriented (or spin) linear bundle isomorphism to the sphere bundle.

The following classical lemma reduces the classification of embeddings to the
relative classification of manifolds (cf. [Sk10, Lemma 1.3]).


\smallskip
{\bf Lemma \Is.}
{\it For a closed connected manifold $N$ two embeddings $f_0,f_1:N\to\R^m$ are
isotopic if and only if there is an oriented bundle isomorphism
$\varphi:\partial C_{f_0}\to\partial C_{f_1}$ which extends to a diffeomorphism $C_{f_0}\to C_{f_1}\#\Sigma$ for some homotopy $n$-sphere $\Sigma$.}

\smallskip
{\it Proof.}
The `only if' part is obvious, so let us prove the `if' part.
The bundle isomorphism $\varphi$ also extends to an
orientation-preserving diffeomorphism
$S^m-\Int C_{f_0}\to S^m-\Int C_{f_1}$.
Therefore $\Sigma\cong S^m\#\Sigma\cong S^m$.
So $\varphi$ extends to an orientation-preserving diffeomorphism
$C_{f_0}\cong C_{f_1}$.
Since any orientation-preserving diffeomorphism of $\R^m$ is isotopic to
the identity, it follows that $f_0$ and $f_1$ are isotopic.
\qed

\smallskip
{\bf Remark \ra.}
{\it Lemma \is has been used to obtain embedding theorems in terms of Poincar\'e
embeddings [Wa70].
But `these theorems reduce geometric problems to algebraic problems which are
even harder to solve' [Wa70].
One of the main problems is that in general (i.e. not in simpler cases like
that of [Sk10, the Effectiveness Theorem]) it is hard to
work with the homotopy type of the pair $(C_f,\partial C_f)$ (which is
sometimes unknown even when the classification of embeddings is known).}

\smallskip
The main idea of our proof
is to apply the
modification of surgery [Kr99] which allows one to classify $m$-manifolds using
their homotopy type just below dimension $m/2$.
\footnote{The realization of this idea is close to, but different from the
realization of [Sk10].
Here we use $B\text{Spin}\times\C P^\infty$-surgery while in [Sk10]
$BO\left<5\right>\times\C P^\infty$-surgery is used.}
Applying modified surgery we prove a diffeomorphism criterion for certain
7-manifolds with boundary: the Almost Diffeomorphism Theorem \hs (cf. the
Diffeomorphism Theorem \Dif) which is a
new, non-trivial version of [KS91, Theorem 3.1] and of [Kr99, Theorem 6] for
7-manifolds $M$ with
non-empty boundary and infinite $H_4(M)$.

\bigskip
{\bf Preparatory results.}

In order to let the reader understand the main ideas before going into details,
we sometimes apply a result before proving it.
In such cases the proof is given in \S3 (except for the proof of `if part' of
the Almost Diffeomorphism Theorem \hs which is given in \S4).

\smallskip
{\bf Remark \rb.} For some readers it would be more convenient to replace
homology by cohomology using Poincar\'e
duality (these readers would
have to pass back to homology at the decisive step of the proof because in
geometric situations like in this paper cup-products are anyway calculated
by passing to cap-products).
For some readers it would be more convenient to replace for a manifold $Q$ a
homology class $z\in H_{n-2}(Q,\partial Q)$ by a homotopy class of a map
$Q\to\C P^\infty$
(then sewing two maps would be a bit more technical) and a spin structure on
$Q$ by a map $Q\to B\text{Spin}$.
See more in [Sk10, Remark 2.3].

\smallskip
Recall that unless otherwise stated:

$\bullet$ {\it $N$ is a closed connected orientable 4-manifold and $f:N\to\R^7$ is an
embedding}.

\smallskip
{\bf Lemma \Nb.} {\it The normal bundle of $f$, $\nu_f$, does not depend on $f$.}

\smallskip
{\it Proof.}
The lemma follows because $\nu=\nu_f$ is completely defined by its
characterictic classes $e$, $w_2$ and $p_1$ [DW59].
We have $e(\nu)=0$, $w_2(\nu)=w_2(N)$ and $p_1(\nu)=-p_1(N)$ by the Wu formulas
because $H_4(N)$ has no torsion.
\qed

\smallskip
Take two embeddings $f_0,f_1:N\to S^7$.
By Lemma \nb there is a bundle isomorphism
$\varphi:\partial C_{f_0}\to\partial C_{f_1}$.
Since $H_1(N)=0$, we have $H^1(\partial C_{f_0})=0$, so $\varphi$
maps the spin structure on $\partial C_{f_1}$ coming from $C_{f_1}\subset S^7$
to the spin structure on $\partial C_{f_0}$ coming from $C_{f_0}\subset S^7$.

By Lemma \is the embeddings $f_0$ and $f_1$ are isotopic if and only if
there is an extension $\overline\varphi:C_{f_0}\to C_{f_1}\#\Sigma$.
Such an extension $\overline\varphi$ sends the generator
$A_{f_0}\in H_5(C_{f_0},\partial)$ to the generator
$A_{f_1}\in H_5(C_{f_1},\partial)$.
Hence $\varphi_*\partial A_{f_0}=\partial A_{f_1}$.

\smallskip
{\bf Agreement Lemma \Agr.} {\it Suppose that $H_1(N)$ has no 2-torsion,
\footnote{We conjecture that this assumption is superfluous when $\varphi$ is a
{\it spin} bundle isomorphism. We conjecture that the converse of the Agreement Lemma \agr holds.}
$f_0,f_1:N\to S^7$ are embeddings and
$\varphi:\partial C_{f_0}\to\partial C_{f_1}$ is an orientation preserving
bundle isomorphism.
We have $\varphi_*\partial A_{f_0}=\partial A_{f_1}$ if
$\varkappa (f_0)=\varkappa (f_1)$.}

\smallskip
Now suppose that $\varkappa (f_0)=\varkappa (f_1)$.
{\it There is a spin bordism between $(C_{f_0},A_{f_0})$ and
$(C_{f_1},A_{f_1})$ relative to the boundaries identified by $\varphi$}
(because by Remark \Rb the obstruction to the existence of such a cobordism
assumes values in $\Omega_7^{Spin}(\C P^\infty)=0$ [KS91, Lemma 6.1]).
It remains to replace the bordism by an $h$-cobordism.
This problem is addressed by modified surgery.
In [Kr99] a surgery obstruction is defined and proved to be complete.
We prove that in our situation the surgery obstruction assumes values in certain Witt group isomorphic to $\Z^4$, i.e. there are four integer-valued
surgery obstructions $\sigma(W)$, $p_1(W)\cdot p_1(W)$, $z^2\cdot z^2$,
$z^2\cdot p_1(W)$ (where $\cdot$ is defined in the Bordism Theorem \Bd).
The heart of our argument is to analyze the dependence of the four surgery
obstructions on choices of the bordism $(W,z)$, homotopy
sphere $\Sigma$ and the bundle isomorphism $\varphi$.
We call the resulting obstruction the Kreck invariant.

\bigskip
{\bf The definition of the Kreck invariant.}


For any manifold $Q$ we abbreviate $H_i(Q,\partial Q)$ to $H_i(Q,\partial)$
and denote Poincar\'e-Lefschetz duality by
$$PD:H^i(Q)\to H_{q-i}(Q,\partial)\quad\text{and}
\quad PD:H_i(Q)\to H^{q-i}(Q,\partial).$$
Recall that for an abelian group $G$ the divisibility $d(0)$ of zero is zero
and the divisibility
$$d(x)\quad\text{of}\quad x\in G-\{0\}\quad\text{is}\quad
\max\{k\in\Z\ | \ \text{there is }x_1\in G:\ x=kx_1\}.$$
A sentence involving $k$ holds for each $k=0,1$.

Take the generator $p\in H^4(B\text{Spin})\cong\Z$ such that  $p = 2 p_1$ where $p_1 \in H^4(B\text{Spin})$ is the pull back of the universal first Pontryagin class in $H^4(BSO)$ (see the proof of Lemma \di in \S3.)
For a compact spin $n$-manifold $W$ take the map $\overline\nu:W\to B\text{Spin}$
corresponding to the given spin structure on $W$ and define
$p_W:=PD\overline\nu^*p\in H_{n-4}(W,\partial)$.

A set $X=(C_0,C_1,A_0,A_1,\varphi)$ consisting of compact connected spin
7-manifolds $C_0$ and $C_1$, generators $A_k\in H_5(C_k,\partial)\cong\Z$ and a spin
diffeomorphism $\varphi:\partial C_0\to\partial C_1$ is called {\it admissible} if
$$\partial A_1=\varphi_*\partial A_0,\quad H_3(\partial C_0)=H_5(\partial C_0)=0,
\quad p_{C_0}=p_{C_1}=0\quad\text{and}\quad d(A_0^2)=d(A_1^2).$$
\quad
According to our strategy we shall define the obstruction $\eta_X$ to extending
$\varphi$ to a diffeomorphism carrying $A_0$ to $A_1$.
\footnote{
A more general situation makes things simpler, but a reader who does not wish
to keep in mind the properties of $C_k,A_k,\varphi$ may assume that
$C_k=C_{f_k}$, $A_k=A_{f_k}$ and $\varphi$ is any spin bundle isomorphism.}

Denote $M_\varphi:=C_0\cup_\varphi(-C_1)$.
For $y\in H_5(M_\varphi)$ and an orientable $n$-submanifold
$C\subset M_\varphi$ we denote
\footnote{
If $y$ is represented by a closed oriented 6-submanifold
$Y\subset M_\varphi$ transverse to $C$, then $y\cap C$ is represented by
$Y\cap C$.
If $C=C_0$, then $y\cap C_0$ is the image of $y$ under the composition of the
homomorphisms $H_n(M_\varphi)\to H_n(M_\varphi,C_1)\to H_n(C_0,\partial)$.
}
$$y\cap C:=PD[(PDy)|_C]\in H_{n-2}(C,\partial).$$

{\bf Null-bordism Lemma \Co.}
{\it Each admissible set has a} null-bordism, {\it i.e. a compact connected
spin 8-manifold $W$ and $z\in H_6(W,\partial)$ such that
$\partial W\underset{Spin}\to= M_\varphi$ and $(\partial z)\cap C_k=A_k$.
Moreover, $\partial z\in H_5(M_\varphi)$ is uniquely defined.}

\smallskip
{\it Proof.}
Look at the segment of (the Poincar\'e-Lefschetz dual to) the
Mayer-Vietoris sequence:
$$H_5(\partial C_0)\to H_5(M_\varphi)\overset{\Psi_1\oplus\Psi_2}\to
\to H_5(C_0,\partial)\oplus H_5(C_1,\partial)\overset{\partial_1-\partial_2}\to
\to H_4(\partial C_0).$$
Here the unmarked arrow is induced by inclusion and $\Psi_kx:=x\cap C_k$.

Since $\partial A_1=\varphi_*\partial A_0$, there is
$A\in H_5(M_\varphi)$ such that $A\cap C_k=A_k$.
Since $H_5(\partial C_0)=0$, such a class $A$ is unique.

Since $\Omega_7^{Spin}(\C P^\infty)=0$ [KS91, Lemma 6.1],
there is a compact spin 8-manifold $W$ and a class $z\in H_6(W,\partial)$
such that $\partial W\underset{Spin}\to=M_\varphi$ and $\partial z=A$.
\qed

\smallskip
Consider the following fragment of the exact sequence of pair
(with any coefficients):
$$H_4(\partial W)\overset{i_W}\to\to H_4(W)\overset{j_W}\to
\to H_4(W,\partial)\overset{\partial_W}\to\to H_3(\partial W).$$
Denote by $\rho_m$ the reduction modulo $m$.

\smallskip
{\bf Definition: the Kreck obstruction $\eta_{W,z}$.}
Take a null-bordism $(W,z)$ of an admissible set $X$.
Denote $d:=d(\partial_Wz^2)$.
Then there is $\overline{z^2}\in H_4(W;\Z_d)$ such that
$j_W\overline{z^2}=\rho_dz^2$.
Define
$$\eta_{W, z}:=\overline{z^2}\cap\rho_d(z^2-p_W)\in\Z_d.$$
\quad
{\it The proof of the independence of $\eta_{W, z}$ of the choice of
$\overline{z^2}$.}
We have $\overline{z^2}-\overline{z^2}'=i_Wa$ for some
$a\in H_4(\partial W;\Z_d)$.
By Lemma \bo below there is $\overline{p_W}\in H_4(W)$ such that
$j_W\overline{p_W}=p_W$.
Then
$$\eta_{W, z}(\overline{z^2})-\eta_{W,z}(\overline{z^2}')=
i_Wa\cap(z^2-p_W)=
i_Wa\cap(\overline{z^2}-\rho_d\overline{p_W})=0.\quad\qed$$

{\bf Lemma \Bo.}
{\it If $(W,z)$ is a null-bordism of an admissible set $X$, then
$\partial_W p_W=0$ and $d(A_0^2)=d(\partial_W z^2)$.}

\smallskip
{\it Proof.}
Consider the segment of the Mayer-Vietoris sequence
$$H_3(\partial C_0)\to H_3(\partial W)
\to H_3(C_0,\partial)\oplus H_3(C_1,\partial)\to H_2(\partial C_0).$$
Since $(\partial_Wp_W)\cap C_k=p_{C_k}=0$ and $H_3(\partial C_0)=0$, we
have $\partial_W p_W=0$.

We have $(\partial_W z^2)\cap C_k=(\partial(z\cap C_k))^2=A_k^2$.
Hence $d(A_k^2)$ is divisible by $d(\partial_Wz^2)$, and
in the above segment of the Mayer-Vietoris sequence $\partial_Wz^2$ is mapped
to $A_0^2\oplus A_1^2$.
If $A_0^2$ is divisible by an integer $d$, then so is $A_1^2$ as well.
Since $H_3(\partial C_0)=0$, we obtain that $\partial_W z^2$ is divisible by
$d(A_0^2)$.
This proves $d(A_0^2)=d(\partial_W z^2)$.
\qed


\smallskip
For an admissible set $X$ by Lemma \bo we can define
$$\eta_X:=\rho_{GCD(A_0^2,24)}\eta_{W,z}\in\Z_{GCD(A_0^2,24)}.$$
\quad
{\it The proof of the independence of $\eta_X$ on the choice of $(W,z)$.}
By the Null-bordism Lemma \co the class $\partial z$ is unique.
The independence of the choice of $(W,z)$ within a spin cobordism class relative
to the boundary is standard (because $p_W$ is a `spin characteristic class').
A change of the spin bordism class  of $W$ (relative to
$\partial W=M_\varphi$) changes $\eta_{W,z}$ by adding
$v^2(v^2-p_1(V))$, where $V$ is some closed spin 8-manifold and
$v\in H_6(V)$.
This is divisible by 24 by the smooth spin case of [KS91, Proposition 2.5].
\qed

\smallskip
{\bf Definition: the Kreck invariant $\eta_u$.}
Assume that $H_1(N)=0$.
Take two embeddings $f_0,f_1:N\to S^7$ such that $\varkappa (f_0)=\varkappa (f_1)=u$.
By Lemma \nb there is a bundle isomorphism
$\varphi:\partial C_{f_0}\to\partial C_{f_1}$.
The different possible spin structures on $\partial C_{f_0}$
are in bijective correspondence with
$H_5(\partial C_{f_0};\Z_2)=H^1(\partial C_{f_0};\Z_2)=0$,
so we may assume that $\varphi$ is spin.
By the Alexander duality, the Agreement Lemma \agr and the fact that $C_k$ are
parallelizable, the set
$X=(C_{f_0}, C_{f_1},A_{f_0},A_{f_1},\varphi)$ is admissible.
Define
$$\eta_u(f_0,f_1):=\eta_X\in\Z_{GCD(u,24)}.$$

This is well-defined because $u=AD(PD(A_{f_0}^2))$ and by the Framing Theorem \Ind($\eta$) below.

For $u\in H_2(N)$ fix an embedding $f_0:N\to\R^7$ such that $\varkappa (f_0)=u$ and
define $\eta_u(f):=\eta_u(f,f_0)$.
(We write $\eta_u(f)$ not $\eta_{f_0}(f)$ for simplicity.)
\footnote{
In general $\eta_u$ depends on the choice of an orientation on $N$,
but $\Emb^7(N)$ by definition does not.}

\bigskip
{\bf The outline of the proof.}


\smallskip
{\it Definition of the framing invariant $\eta'_X$.}
Take an admissible set $X=(C_0,C_1,A_0,A_1,\varphi)$ such that
$A_0^2$ and $A_1^2$ are divisible by 2.
Define $\overline{z^2}\in H_4(W;\Z_2)$ analogously to
$\overline{z^2}\in H_4(W;\Z_d)$ in the definition of $\eta_{(W, z)}$.
Define
\footnote{This is independent of the choice of $(W,z)$ analogously to $\eta_X$ using
the smooth spin case of [KS91, Proposition 2.5] (because $12S_3-48S_2=6z^4$
is divisible by 12, so $z^4$ is divisible by 2 for closed manifolds).}
$$\eta'_X:=\overline{z^2}\cap\rho_2z^2\in\Z_2.$$

{\bf Almost Diffeomorphism Theorem \Hs.}
{\it Let $X=(C_0,C_1,A_0,A_1,\varphi)$ be an admissible set such that
$\pi_1(C_k)=H_3(C_k)=H_4(C_k,\partial)=0$.
For some homotopy 7-sphere $\Sigma$ there is a diffeomorphism
$\overline\varphi:C_0\to C_1\#\Sigma$ extending $\varphi$ and such that $\overline\varphi_*A_0=A_1$ if and only if
$$\eta_X=0\quad\text{and, for $A_0^2$ divisible by 2,}\quad \eta'_X=0.$$}
\quad
The `only if' part is simple.
(Indeed, take a 3-connected almost parallelizable 8-manifold $V$ such that
$\partial V=-\Sigma$.
Define $W:=C_0\times I\sharp V$.
Then $\partial W=C_0\cup(-C_0)\#\partial V\cong C_0\cup_\varphi(-C_1)$.
Define $z:=A_0\times I\sharp0$.
Then $(\partial z)\cap C_1=A_1$ because
$\overline\varphi_*A_0=A_1$.
We have $p_W=p_{C_0}\times I+p_V=0$ and $z^4=A_0^4\times I+0=0$.
Thus $\eta_X=0$ and, for $A_0^2$ divisible by 2, $\eta_X'=0$.)
This part is not used in the proof of the Classification Theorem \Ct.


\smallskip
{\bf Framing Theorem \Ind.}
{\it Let $X=(C_0,C_1,A_0,A_1,\varphi)$ be an admissible set
such that $\partial C_k$ is an $S^2$-bundle over a closed
4-manifold $N$ and $\varphi$ is a spin bundle isomorphism.
Then}

($\eta$) {\it $\eta_X$ is independent of the choice of bundle isomorphism $\varphi$ (the choice preserving $C_k$,
$A_k$ and admissibility).
\footnote{The change of $\varphi$ is only possible together with certain
changes of $(W,z)$.}}

($\varphi$) {\it If $A_0^2$ is divisible by 2, then we can change bundle isomorphism $\varphi$
(change preserving $C_k$, $A_k$ and admissibility) so as to obtain $\eta'_X=0$.}


\smallskip
{\bf Transitivity Lemma \Ad.}
{\it If $f,f_1,f_2:N\to\R^7$ are embeddings with the same value of the
Bo\'echat-Haefliger invariant, $u$, then
$\eta_u(f,f_1)+\eta_u(f_1,f_2)=\eta_u(f,f_2)$.}

\smallskip
{\it Proof of the injectivity of $\eta_u$.}
By the Transitivity Lemma \ad it suffices to prove the following:

$\bullet$ {\it if $\varkappa (f)=\varkappa (f')$ and $\eta_{\varkappa (f)}(f,f')=0$, then $f$ is isotopic to
$f'$.}

In order to prove this assertion construct an admissible set $X$ as in the
definition of the Kreck invariant $\eta_u(f,f')$.
Since $\eta_u(f,f')=0$, we have $\eta_X=0$.

If $A_f^2$ is divisible by 2, then by the Framing Theorem \Ind($\varphi$) we can
change $\varphi$ so as to obtain $\eta'_X=0$.
By the Framing Theorem \Ind($\eta$) $\eta_X$ will be preserved.

Therefore by the Almost Diffeomorphism Theorem \hs $\varphi$ extends to a
diffeomorphism $C_f\to C_{f'}\#\Sigma$ for a certain homotopy 7-sphere $\Sigma$.
Hence $f$ is isotopic to $f'$ by Lemma \Is.
\qed

\smallskip
{\it The description of $\im\eta_u$} holds by the second equality of the
Addendum \adm and the following two partially known results proved in \S3.

\smallskip
{\bf Lemma \Di.}
{\it Let $W$ be a compact spin $8$-manifold.
Then

(a) $2p_W=PDp_1(W)$.

(b) $p_W\cap x-x\cap x$ is divisible by 2 for each $x\in H_4(W)$.}

\smallskip
{\bf Realization Theorem \Re.}
{\it There is an embedding $g_1:S^4\to S^7$ such that $\eta_0(g_1)=2$.}

\smallskip
The Realization Theorem \re holds by the injectivity of $\eta_0$ (proved above)
because there exist 12 pairwise non-isotopic embeddings $S^4\to S^7$ [Ha66].
We present an alternative direct proof in \S3.

In what follows please note that Sections \S3 and \S4 depend on \S2 but are independent of each other.

\head 3. Further details for the proof \endhead

{\bf Proof of the Agreement Lemma \Agr.}

For a map $\xi:P\to Q$ between a $p$-manifold and  a $q$-manifold denote the `preimage'
homomorphism by
$$\xi^!:=PD\circ\xi^*\circ PD:H_i(Q,\partial)\to H_{p-q+i}(P,\partial).$$
Let $f_0:N\to S^7$ be an embedding.
In this subsection we omit the subscript $f_0$ from $\nu_{f_0}$, $C_{f_0}$, $A_{f_0}$ etc.

Let $N_0:=\Cl(N-B^4)$, where $B^4$ is a closed 4-ball in $N$.
Let $\zeta:N_0\to \nu^{-1}N_0$ be a section of the normal bundle
$\nu^{-1}N_0\to N_0$.
(This exists because $e(\nu)=0$.)
Consider the following diagram.
$$\minCDarrowwidth{0pt}\CD
H_4(N_0,\partial) @>> \zeta_* > H_4(\nu^{-1}N_0,\partial) @<< e <
H_4(\partial C,\nu^{-1}B^4) @<< j <  H_4(\partial C) @>> i > H_4(C)\endCD.$$
Here $j$ is the isomorphism from the exact sequence of a pair, $e$ is the
excision isomorphism and $i$ is induced by the  inclusion.
For $k\ne0$ we identify $H_k(N)$ with $H_k(N_0,\partial)$
by the composition
$H_k(N)\overset{j_N}\to\to H_k(N,B^4)\overset{e_N}\to\to H_k(N_0,\partial)$
of the isomorphism from the exact sequence of the pair $(N, B^4)$ and the excision
isomorphism.

Consider the following fragment of the Gysin sequence for the bundle
$\nu$ having trivial Euler class:
$$0\to H_2(N)\overset{\nu^!}\to\to H_4(\partial C)\overset{\nu_*}\to
\to H_4(N)\to0.$$
We see that the map
$$\nu_*\oplus\zeta^!ej:H_4(\partial C)\to H_4(N)\oplus H_2(N)$$
is an isomorphism.
By the definitions of $A$ and $A_{f_1}$ we have
$\nu_*\partial A=[N]=\nu_{f_1*}\partial A_{f_1}$.
So it suffices to prove that
$$(*)\qquad\zeta^!ej\partial A=(\varphi\zeta)^!e_{f_1}j_{f_1}\partial A_{f_1}
\quad\text{for some section}\quad\zeta:N_0\to \nu^{-1}N_0.$$
We shall call a section $\zeta$ {\it weakly unlinked} if $ij^{-1}e^{-1}\zeta_*=0$.

\smallskip
{\bf Section Lemma \Eu.} [Sk10, Section Lemma 2.5.b]
{\it If $\zeta$ is a weakly unlinked section, then
$\varkappa (f)=PDe(\zeta^\perp)=\zeta^!ej\partial A$, where $\zeta^\perp$ is the
oriented $S^1$-bundle that is the orthogonal complement to $\zeta$ in
$\nu|_{N_0}$.}

\smallskip
There exist unlinked sections $\zeta_0$ and $\zeta_1$ for $f_0$ and $f_1$
[HH63, 4.3, BH70, Proposition 1.3, Sk08', the Unlinked Section Lemma (a)]
(because by [Sk10, Remark 2.4 and footnote 14] our definition of the weakly
unlinked section is equivalent to the original definition [BH70]).
By the Section Lemma \eu (*) is implied by $\varphi\zeta_0=\zeta_1$.
For sections
$$\xi,\eta:N_0\to\partial C_{f_1}\quad\text{we have}\quad
PDe(\xi^\perp)-PDe(\eta^\perp)=\pm2d(\xi,\eta),$$
where $d(\xi,\eta)\in H_2(N)$ is the difference element [BH70, Lemme 1.7].
Since $H_2(N)$ has no 2-torsion, $d(\varphi\zeta_0,\zeta_1)=0$ follows from
$$PD e((\varphi\zeta_0)^\perp)=PDe(\zeta_0^\perp)=\varkappa (f_0)=\varkappa (f_1)=PDe(\zeta_1^\perp),$$
where the second and the last equalities holds by the first equality of the Section Lemma \Eu.
\qed


\bigskip
{\bf Proof of the Framing Theorem \Ind.}

\smallskip
{\bf Lemma \Ide.}
{\it Define $i:S^1=SU_1\to SU_3$ by $i(z)=\diag(z,\overline z,1)$.
Then the homogeneous space $SU_3/i(S^1)$ is the total space of the
non-trivial $S^2$-bundle over $S^5$ (i.e. the bundle corresponding to the
non-trivial element of $\pi_4(SO_3)\cong\Z_2$).}

\smallskip
{\it Proof.} Since $i(S^1)\subset SU_2$, the standard bundle
$SU_2\to SU_3\to S^5$ gives a bundle
$$(*)\qquad S^2\cong SU_2/i(S^1)\to SU_3/i(S^1)\to S^5.$$
Here the diffeomorphism is given by a free action of $SU_2$ on $\C P^1=S^2$
whose stabilizer subgroup is $i(S^1)$.

(In order to define such an action, identify $SU_2$ with the group of unit
length quaternions.
Define the Hopf map
$$h:SU_2\to\C P^1\quad\text{by}\quad h(z+jw):=(z:w)\quad
\text{for}\quad z,w\in\C\quad\text{and}\quad |z|^2+|w|^2=1.$$
The required action is well-defined by $uh(v):=h(uv)$.
The action of $SU_2$ on $\C^2=\Ha$ is given by
$(z+jw)(p+jq)=zp+\overline wq+j(wp+\overline zq)$.
Hence $z+jw$ corresponds to the
matrix $\left(\matrix z &w\\ -\overline w & \overline z\endmatrix\right)$.
Thus the stabilizer subgroup is $\{z+j0\ |\ z\in\C\}=i(S^1)$.)

Since $\pi_4(SU_3)=0$
(by $\pi_4(SU_3)\cong\pi_4(SU)$ and the Bott periodicity), we have
\linebreak
$\pi_4(SU_3/i(S^1))=0\ne\Z_2\cong\pi_4(S^2\times S^5)$.
Hence $SU_3/i(S^1)\not\cong S^2\times S^5$.
Therefore the bundle (*) is non-trivial.
\qed
\footnote{{\it An alternative proof of the non-triviality of the bundle (*).}
If (*) is trivial, then there is a bundle $S^1\to SU_3\to S^2\times S^5$
whose first Chern class is a generator of $H^2(S^2\times S^5)\cong\Z$.
Then $SU_3\cong S^3\times S^5$ which is a contradiction because
$\pi_4(SU_3)=0\ne\Z_2\cong\pi_4(S^3\times S^5)$.}

\smallskip
{\it Proof of the Framing Theorem \Ind.}
Take a closed 4-ball $B^4\subset N$.
Since $\partial C_0$ is an $S^2$-bundle over a closed 4-manifold $N$ and $H_3(\partial C_0)=0$, we have $H_1(N)=H_3(N)=0$.
This and $\pi_2(SO_3)=0$ imply that the spin bundle isomorphism $\varphi$ is uniquely defined over $\Cl(N-B^4)$.
If we change $\varphi$ over $B^4$, then analogously to [Sk08', proof of the Independence Lemma] and by Lemma \ide the pair $(M_\varphi,A_0\cup_\varphi A_1)$ would change
by connected sum over $S^2$ with $(SU_3/i(S^1),A)$, where $A\in H_5(SU_3/i(S^1))\cong\Z$.
It suffices to consider the case when $A$ is a generator.

We have that $SU_3/i(S^1)$ is $N_{1,-1}$ defined in [KS91, \S1]; the assumption
$k+l\ne0$ is not used for the definition (but it is required for the positive
curvature properties Kreck and Stolz consider).
By [KS91, Proposition 2.2]
$(SU_3/i(S^1),A)\underset{Spin}\to=\partial(W,z)$ for some spin
8-manifold $W$ and $z\in H_6(W,\partial)$.
By Lemma \ide $H_3(\partial W)=H_4(\partial W)=0$.
Hence we may identify $z^2$ and $p_W$ with elements of $H_4(W)$ (these elements are
denoted by the same letters).
In [KS91, proof of Lemma 4.4] the assumption $k+l\ne0$ was not used.
\footnote{There is a typographical error in the expression for $s_3$ which
should read $s_3(N_{k,l})=(-4P+NS)/6N$ and in the expression for $P$ where
$-6m^2n^2$ should read $-6lm^2n^2$; we do not use these corrections.}
So by [KS91, (2.4), Lemma 4.4 and the bottom of p. 475] with
$$k=m=1,\quad l=-1,\quad n=0\quad\text{we have}\quad z^4=-1\quad\text{and}
\quad N=P=S=1,\quad$$
$$\text{so}
\quad -2z^2p_W+2z^4=48s_2(N_{1,-1})=2(-P+NS)/N=0.$$
Thus any change of $\varphi$
preserves $\eta_X$ and, for $A_0^2$ divisible by 2, there is a change of $\varphi$ that changes $\eta'_X$ by 1.
\qed

\bigskip
{\bf Proof of the the Transitivity Lemma \Ad.}

Assume that $(W_k,z_k)$ is a null-bordism of the admissible set
$(C_f, C_{f_k},A_f,A_{f_k},\varphi_k)$.

Take $\varphi:=\varphi_2\varphi_1^{-1}$.
Then $X=(C_{f_1},C_{f_2},A_{f_1},A_{f_2},\varphi)$ is admissible.

Take $W:=W_2\cup_{C_f}(-W_1)$.
From the Mayer-Vietoris sequence
$$H_6(C_f)\to H_6(W,\partial)\overset\Psi\to
\to H_6(W_1,\partial)\oplus H_6(W_2,\partial)\to H_5(C_f)$$
we see that $\Psi$ is an isomorphism.
Take $z:=\Psi^{-1}(z_1\oplus z_2)$.
Then $(W,z)$ is a null-bordism of $X$.

Denote $d:=d(A_0^2)=d(A_1^2)$.
Consider the maps
$$(\cdot\cap W_1)\oplus(\cdot\cap W_2):H_4(W,\partial)
\to H_4(W_1,\partial)\oplus H_4(W_2,\partial)\quad\text{and}$$
$$i_1\oplus i_2:H_4(W_1;\Z_d)\oplus H_4(W_2;\Z_d)\to H_4(W;\Z_d).$$
Clearly, $p_{W_k}=p_W\cap W_k$ and $z^2_k=z^2\cap W_k$.
Take $\overline{z^2}:=i_1\overline{z_1^2}\oplus i_2\overline{z_2^2}$.
Since
$$(i_1x_1\oplus i_2x_2)\cap y=-x_1\cap(y\cap W_1)+x_2\cap(y\cap W_2)\quad
\text{we have}\quad\eta_{W, z}=-\eta_{W_1,z_1}+\eta_{W_2,z_2}.$$
Hence $\eta_u(f_1,f_2)=-\eta_u(f,f_1)+\eta_u(f,f_2)$.
\qed

\bigskip
{\bf Proof of the second equality of the Addendum \Adm.}

It suffices to prove that
$$\eta_u(f\#g,f_0\#g_0)=\eta_u(f,f_0)+\eta_0(g,g_0),$$ where
$f_0:N\to S^7$ is any embedding, $u=\varkappa (f_0)=\varkappa (f)$ and $g_0:S^4\to\R^7$ is the standard embedding.
By the Null-Bordism Lemma \nb there is a null-bordism
$(W_f,z_f)$ of an admissible set
$(C_f, C_{f_0},A_f,A_{f_0},\varphi_f)$.
Analogous assertion holds with $f, f_0$ replaced by $g, g_0$.

Since $H_1(N)=H_3(N)=0$ and $\pi_2(SO_3)=0$, we may assume that $\varphi_f$ is the identity outside $B^4\subset N$ and
that $\nu_f=\nu_{f\#g}$ outside $B^4\subset N$.
Then take any spin bundle isomorphism
$\varphi:\partial C_{f\#g}\to\partial C_{f_0\#g_0}$
that is the identity outside $B^4$.

Identify $B^4\times S^2$ and $\nu_f^{-1}B^4\subset\partial C_f$ and do the same for $f$ replaced by $f_0,g$ or $g_0$.
We have
$$C_{f\#g}=C_f\cup_{B^4\times S^2}C_g\quad\text{and}
\quad C_{f_0\#g_0}=C_{f_0}\cup_{B^4\times S^2}C_{g_0}.$$
Then $(C_{f\#g},C_{f_0\#g_0}, A_f, A_{f_0}, \varphi)$ is an admissible set.

By $B^5=B^5_+\cup_{B^4}B^5_-$ we denote the standard decomposition.
Take an embedding $B^5\times S^2\to\partial W_f=C_f\cup_{\varphi_f}C_{f_0}$
whose image intersects
$$C_f,\quad C_{f_0}\quad\text{and}\quad
\partial C_f\overset{\varphi_f}\to=\partial C_{f_0}\quad\text{by}
\quad B^5_+\times S^2,\quad B^5_-\times S^2\quad\text{and}\quad B^4\times S^2,$$
respectively.
Take the analogous embedding $B^5\times S^2\to\partial W_g$.
Then take
$$W:=W_f\cup_{B^5\times S^2}W_g.$$
Consider the Mayer-Vietoris sequence:
$$H_6(B^5\times S^2)\to H_6(W,\partial)
\to H_6(W_f,\partial)\oplus H_6(W_g,\partial)\to H_5(B^5\times S^2,\partial).$$
Identify $\partial W$ and $C_{f\#g}\cup_\varphi C_{f_0\#g_0}$ by the easily
constructed homeomorphism.
We have $\partial A_f\cap (B^4\times S^2)=[B^4\times x]
\in H_4(B^4\times S^2,\partial)$, and the same for $f$ replaced by $f_0,g$ or $g_0$.
Hence
$$\partial z_f\cap (B^5\times S^2)=\partial z_g\cap (B^5\times S^2)=
[B^5\times x]\in H_5(B^5\times S^2,\partial).$$
Therefore there is a unique $z\in H_6(W,\partial)$ such that $(W,z)$ is a
null-bordism of
\linebreak
$(C_{f\#g},C_{f_0\#g_0}, A_f, A_{f_0}, \varphi)$.

Since $H_i(B^5\times S^2)=H_i(B^5\times S^2,\partial)=0$ for $i=3,4$,
by the homology exact sequence of the pair and the Mayer-Vietoris sequence we have
isomorphisms $\Psi$ and $\Psi_\partial$; these isomorphisms are also isomorphisms of the respective intersection forms and fit into the following commutative diagram:
$$\minCDarrowwidth{5pt} \CD  H_4(W) @<< \Psi\cong < H_4(W_f)\oplus H_4(W_g) \\
  @VV j V  @VV j_f\oplus j_g V \\
H_4(W,\partial) @>> \Psi_\partial\cong
> H_4(W_f,\partial)\oplus H_4(W_g,\partial).
\endCD$$
Clearly, $\Psi_\partial z^2=z^2_f\oplus z^2_g$ and
$\Psi_\partial p_W=p_{W_f}\oplus p_{W_g}$.
So we can take
$\overline{z^2}:=\Psi_d(\overline{z^2_f}\oplus\overline{z^2_g})$, where
$\Psi_d$ denotes the isomorphism analogous to $\Psi$ with coefficients $\Z_d$.
Then clearly $\eta_{W, z}=\eta_{W_f,z_f}+\eta_{W_g,z_g}$.
This implies the required statement.
\qed
\footnote{We conjecture that
$\eta_{u_1\oplus u_2}(f_1\#f_2,f_1'\#f_2')= \rho_{GCD(u_1,u_2,24)}\eta_{u_1}(f_1,f_1')+\rho_{GCD(u_1,u_2,24)}\eta_{u_2}(f_2,f_2')$,
where $f_k,f_k':N_k\to\R^7$ are embeddings such that $\varkappa (f_k)=\varkappa (f_k')=u_k$.}


\bigskip
{\bf Proof of Lemma \Di.}

\smallskip
{\it Proof of part} (a).
Consider the fibration $\R P^\infty\to B\text{Spin}\to BSO$.
The 4-line of the cohomology Leray-Serre spectral sequence of this fibration is
the same at the $E_2$ term and at the $E_\infty$ term.
The 4-line has $\Z=H^4(BSO)$ in the $(4,0)$ position and also a
$\Z_2=H^2(BSO;\Z_2)$ in the $(2, 2)$ position.
Therefore $H^4(BSO)$ maps into $H^4(B\text{Spin})$ as a subgroup of index 2.
Hence the pullback $p_1\in H^4(B\text{Spin})$ of the universal first Pontryagin
class in $H^4(BSO)$ equals $2p$. (This is also proved in [KS91, proof of Lemma 6.5].)
Then $2p_W=PD\overline\nu^*p_1=PDp_1(W)$.
\qed

\smallskip
{\it Proof of (b).}
Let $w_4\in H^4(B\text{Spin};\Z_2)$ be the pullback of the universal 4-th
Stiefel-Whitney class in $H^4(BSO;\Z_2)$.
Since $w_4$ generates $H^4(B\text{Spin};\Z_2)$ and the mod 2 reduction
$\rho_2:H^4(B\text{Spin})\to H^4(B\text{Spin};\Z_2)$ is onto, we have $\rho_2(p)=w_4$.
Also $w_4(W)=\overline\nu^*w_4$.
Hence $\rho_2(p_W)=PDw_4(W)$.
Let us prove that this implies (b).

If $W$ is closed, then part (b) follows
because $w_4(W)=v_4(W)+\Sq^1 v_3(W)=v_4(W)$.
Here the first equality holds by the Wu formula and the second because
$\Sq^1 v_3(W)=\Sq^1 w_3(W)=0$ since $W$ is spin
(or else because $v_3(W)=w_3(W)=0$ since $W$ is spin and the space $B\text{Spin}$ is
3-connected).

If $W$ has a non-empty boundary, then let $Y:=W\cup_{\partial W}(-W)$.
Since
$$p_W =p_Y\cap W,\quad\text{we have}\quad
p_W\cap_W x=p_Y\cap_Y i_Yx\underset{\mod2}
\to\equiv i_Yx\cap_Y i_Yx=x\cap_W x,$$
where $i_Y$ is the inclusion-induced map $H_4(W)\to H_4(Y)$.
\qed


\bigskip
{\bf Proof of the Realization Theorem \Re.}

\smallskip
{\it A construction of $g_1:S^4\to S^7$.}
By general position, there is an embedding $\eta'':S^3\to S^2\times D^5$
whose composition with the projection onto $S^2$ is the Hopf map.
\footnote{{\it An explicit construction of $\eta''$.}
Define an embedding $\eta':S^3\to S^2\times D^2$ by
$\eta'(z_1,z_2):=((z_1:z_2),z_1)$.
The composition of $\eta'$ with the projection onto $S^2$ is the Hopf map.
Let $\eta''$ be the composition of $\eta'$ and the standard inclusion
$S^2\times D^2\to S^2\times D^5$.}
Take an embedding $\psi:D^4\to S^2\times D^5$ whose image intersects
$\eta''(S^3)$ transversally at exactly one point of sign $+1$.
Let $\psi':=\psi|_{\partial D^4}$.

Since each embedding $\alpha:S^3\to S^7$ is unknotted, it extends to an embedding
$D^4\to D^8\supset S^7$.
Since $D^4$ is contractible, it has a unique framing.
Therefore there is a unique framing of $\alpha(S^3) \subset S^7$ which extends to a
framing of some extension $D^4\to D^8$.
Define this framing to be the zero framing.
This and the isomorphism $\pi_3(SO_4)\cong\Z\oplus\Z$ [Mi56] give
{\it a 1--1 correspondence between normal framings of embedding $\alpha:S^3\to S^7$
(up to homotopy) and $\Z\oplus\Z$.}

Assume that $S^2\times D^5\subset S^7$ is standardly embedded as a complement
to the tubular neighborhood of the standard $S^4 \subset S^7$.
Take the framing on $\eta''$ corresponding to $(0,0)$ and
the framing on $\psi'$ corresponding to $(1,-1)$.
Let $M$ be the closed 7-manifold obtained from $S^7$ by surgery along
framed embeddings $\psi'$ and $\eta''$.
In the `proof of the Realization Theorem \Re' below we prove that $M\cong S^7$.
Let $g_1$ be the composition of the inclusion $S^4\to M$ and any
diffeomorphism $M\to S^7$.

\smallskip
In this subsection let $i: S^2 \times D^5_{-} \to S^7=\partial D^8$
be the standard embedding.
For a $D^4$-bundle $\widetilde\alpha$ over $S^4$ denote by
$e(\widetilde\alpha),p_1(\widetilde\alpha)\in\Z$ the Euler and the Pontryagin numbers of this bundle (defined using the standard orientation on $S^4$).

\smallskip
{\bf Lemma \Ir.} {\it Let $W$ be the $8$-manifold obtained by
adding $4$-handles to $S^2\times D^6$ via embeddings
$$\alpha_1,\dots,\alpha_n:S^3\times D^4\longrightarrow S^2\times D^5_-\subset
\partial(S^2\times D^6)$$
with disjoint images.
Denote by $[\alpha_1],\dots,[\alpha_n]\in H_4(W)$ the basis corresponding to
the $4$-handles.
Denote by $\widetilde\alpha_m$ the $D^4$-bundle over $S^4$ corresponding to
$\alpha_m$ (i.e. the projection from the 8-manifold $W'_m$ obtained from
$D^8$ by adding a 4-handle along $i\alpha_m$ to the sphere $S^4_m\subset W'_m$ representing $i\alpha_m$).
Then
$$[\alpha_m]\cap[\alpha_l]=
\cases\lk_{S^7}(i\alpha_m,i \alpha_l)&m\ne l\\
e(\widetilde\alpha_m)&m=l\endcases\quad\text{and}\quad
2p_W\cap[\alpha_m]=p_1(\widetilde\alpha_m).$$}

{\it Proof.} Cf. [Sc02].  The equality
$[\alpha_m]\cap[\alpha_l]=\lk_{S^7}(\alpha_m,\alpha_l)$ for $l\ne m$ follows
analogously to [Ma80, 3.2].
For the other equalities we may assume that $m=l=1$, replace $W$ by $W'_1$ and omit subscripts 1.


We have
$[\alpha]\cap[\alpha]=e(\widetilde\alpha)$ because the self-intersection of a
homology class represented by a submanifold equals to the Euler class of the
normal bundle of the submanifold in the manifold (this is easily proved
directly or else deduced from [MS74, Exercise 11-C in p. 134]).

We have
$2p_{W'}\cap[\alpha]=PDp_1(\tau_{W'}|_{S^4})=PDp_1(\widetilde\alpha)$, where
the second equality holds because
$\tau_{W'}|_{S^4}\cong\tau_{S^4}\oplus\nu_{W'}(S^4)$ is stably equivalent to
$\nu_{W'}(S^4)=\widetilde\alpha$ since $S^4$ is stably parallelizable.
\qed

\smallskip
{\it Proof of the Realization Theorem \Re.}
Let $S^2\times\partial D^6=S^2\times D^5_+\cup_{S^2\times S^4}S^2\times D^5_-$
be the standard decomposition corresponding to the standard decomposition
$\partial D^6=D^5_+\cup_{S^4}D^5_-$.
Let $W$ be the 8-manifold obtained from $S^2\times D^6$ by adding $4$-handles
along the framed embeddings $\psi'$ and $\eta''$ into $S^2\times D^5_-$.
Let $C_0:=S^2\times D^5_+\subset\partial W$.
Let $C_1\subset\partial W$ be the 7-manifold obtained from $S^2\times D^5_-$
by surgery along framed embeddings $\psi'$ and $\eta''$ into $S^2\times D^5_-$.
Take the identity diffeomorphism $\varphi:\partial C_0\to\partial C_1$.


Take a basis $x,y$ of $H_4(W)\cong\Z^2$ with $x$ and $y$ corresponding to the
handle attached by $\psi'$ and by $\eta''$, respectively.
By Lemma \ir and [Mi56]
$$x\cap y=1,\quad x\cap x=p_W\cap x=0,\quad y\cap y=1+(-1)=0
\quad\text{and}\quad p_W\cap y=1-(-1)=2.$$
Hence $p_W=2x$.

In this paragraph we prove that {\it $X=(C_0,C_1,A_0,A_1,\varphi)$
is an admissible set and $(W,z_W)$ is a null-bordism of $X$.}
For the $W$ we constructed above the maps of the composition
$H_6(W,\partial)\to H_5(\partial W)\to H_5(C_k,\partial)$,
the boundary map and
the map $x\mapsto x\cap\partial C_k$, are both isomorphisms.
Hence for the generator $z_W\in H_6(W,\partial)$ we have that $\partial z_W$
is a generator of $H_5(\partial W)$ and that $A_k:=\partial z_W\cap C_k$
is a generator of $H_5(C_k,\partial)$.
Clearly, $p_{C_0}=0$.
Since the intersection form $H_4(W)\times H_4(W)\to\Z$ is non-degenerate,
the map $j:H_4(W)\to H_4(W,\partial)$ is an isomorphism.
This and $H_3(W)=0$ imply by the exact sequence of the pair $(W,\partial W)$
that $H_3(\partial W)=0$.
Since the inclusion $H_2(\partial C_0)\to H_2(W)$ is an isomorphism,
using Mayer-Vietoris sequence we obtain that $H_3(C_1,\partial)=0$.
Hence $p_{C_1}\in H_3(C_1,\partial)=0$.


Denote by $W'$ the 8-manifold obtained from $D^8$ by adding 4-handles along
framed embeddings $i\psi'$ and $i\eta''$ into $\partial D^8$.
Recall that $M=\partial W'$ for the 7-manifold $M$ defined in the `construction
of $g_1$'.
Analogously to above there is a basis $x,y$ of $H_4(W')\cong\Z^2$
in which the intersection form of $W'$ has matrix $H_+:=\left(\matrix0&1\\1&0\endmatrix\right)$, and $p_{W'}=2x$.
Then $\sigma(W')=0\underset{\mod28\cdot8}\to\equiv0=p_{W'}\cap p_{W'}$.
Hence $\partial W'\cong S^7$ [EK62, \S6].

We have $z_W^2=y$.

(Indeed, $W\simeq S^2 \cup (e^4_x \cup e^4_y)$, where $\simeq$ means
`homotopy equivalent up to dimension 4'.
Homotopy classes of the attaching maps for $e^4_x$ and for $e^4_y$ equal to
the homotopy classes of $\eta''$ and $\psi'$.
So the attaching maps are homotopic to the Hopf map and trivial map
$S^3\to S^2$, respectively.
It follows that $W \simeq \C P^2 \vee S^4$.
Thus we obtain the cohomology ring of $W$ up to dimension 4.
By duality we obtain the homology groups of $W$ and
relevant intersection products above dimension 3.
Hence $z_W^2\cap x=1$ and $z_W^2\cap y=0$ for a generator $z_W\in H_6(W)$.
By Poincar\'e duality $z_W^2=y$.)

Then $\eta(g_1,g_0)=\eta_{W,z_W}=2$.
\qed

\head 4. Proof of the `if' part of the Almost Diffeomorphism Theorem
\Hs \endhead

\smallskip
{\bf The Kreck Theorem \Kr.}
{\it Let

$\bullet$ $W$ be a compact $4l$-manifold such that $\partial W=C_0\cup C_1$ for
compact $(4l-1)$-manifolds $C_0,C_1\subset\R^{8l}$ with common boundary;

$\bullet$ $p:B\to BO$ be a fibration such that $\pi_i(p)=0$ for $i\ge 2l$ and
$\pi_1(B)=0$;

$\bullet$ $\overline\nu:W\to B$ be a $2l$-connected map such that $p\overline\nu|_{C_k}$
is the classifying map of the normal bundle of $C_k$ and $\overline\nu|_{C_k}$
is $(2l-1)$-connected.

Then $\overline\nu$ is bordant (relative to the boundary) to a product
of $\overline\nu|_{C_0}$ with the interval if
\footnote{The `only if' implication also holds but is not used in this paper
(the same is true for the Bordism Theorem \bd below).}
there is a subgroup $U\subset H_{2l}(W)$ such that

$\bullet$ $U\cap U=0$ and $\overline\nu_*U=0\subset H_{2l}(B)$,

$\bullet$ $j_k|_U$ is an isomorphism onto a direct summand
in $V_k:=H_{2l}(W,C_k)$,
and

$\bullet$
the quotient $j_0U\times V_1/j_1U\to\Z$ of the intersection pairing
$\cap:V_0\times V_1\to\Z$ is unimodular. }

\smallskip
{\it Proof.} Denote $K:=\ker(\overline\nu_*:H_{2l}(W)\to H_{2l}(B))$.
The form $\cap:K\times K\to\Z$ is even because
\footnote{In the situation of the Almost Diffeomorphism Theorem \hs this
form is even by Lemma \Di(b).}
$$x\cap x=\left<v_{2l}(W),x\right>=\left<p^*\overline\nu^*v_{2l},x\right>=
\left<v_{2l},p_*\overline\nu_*x\right>=0\mod2,$$
where $x\in K$
and $v_{2l}\in H^{2l}(BO)$ is the $2l$th
Wu class.
So in [Kr99, p. 725] we can take $\mu(x):=x\cap x/2$ for $x\in K$ (because $2l$
is even).
We have $Wh(\pi_1(B))=0$ and so an isomorphism is a simple isomorphism.
Hence the hypothesis on $U$ implies that $\theta(W,\overline\nu)$ is
`elementary omitting the bases' [Kr99, Definition in p. 730 and the second
remark on p. 732].
\footnote{In [Kr99, Definition on p. 729] $\theta(W,\overline\nu)$ was only
defined for a $q$-connected map $\overline\nu:W\to B$.
(Indeed, on p. 725 in [Kr99]
there is a paragraph beginning {\it ``The objects in $l_{2q}(\pi, \omega)$ are
represented ... "}.
In condition (i) $V_0$ and $V_1$ are based.
This means in particular that they are stably free.   Now for a
bordism $(W,\overline\nu; M_0, M_1)$ we have by definition
$V_0 = H_q(W,M_0)$ and this is only a stably free module if
$\overline\nu : W \to B$ is $q$-connected.)
If $\overline\nu$ is not $q$-connected, then it is bordant to a $q$-connected
map $\overline\nu_1:W_1\to B$ and we can define
$\theta(W,\overline\nu):=\theta(W_1,\overline\nu_1)$.
This is well-defined by [Kr99, the first sentence in p. 730]. }
Thus the result follows by the $h$-cobordism theorem and [Kr99, Theorem 3 and
second remark in p. 732].
\qed

\bigskip
{\bf The Bordism Theorem.}


\smallskip
{\bf Lemma \Ho.}
{\it For $k=0,1$ let $C_k$ be compact connected 7-manifolds such that $H_3(C_k)=0$, let
\ $\varphi:\partial C_0\to\partial C_1$ be a diffeomorphism and let $W$ be a compact
8-manifold such that $\partial W=M_\varphi$.
Denote
$$V_0:=H_4(W,C_0)\quad\text{and let}\quad
j_0:H_4(W)\to V_0$$
be the map from the exact sequence of the pair $(W, C_0)$.  Then there is a well-defined bilinear map
$$\cdot:V_0\times V_0\to\Z\quad\text{given by}
\quad x\cdot x':=j_0^{-1}x\cap x'$$
which is symmetric and unimodular: here $j_0^{-1}x$ denotes any element in $j_0^{-1}x$.}

\smallskip
{\it Proof.}
Since $H_3(C_0)=0$, the map $j_0$ is surjective.

If $y,y'\in j_0^{-1}x$, then we may assume that the support
of $y-y'$ is in $C_0$.
Then $(y-y')\cap x'=(y-y')\cap_{C_0}\partial x'=0$ because $H_3(C_0)=0$.
So $\cdot$ is well-defined.

This form is symmetric because of the symmetry of linking coefficients of
3-cycles in $C_0$.
In order to prove the unimodularity of $\cdot$ take a primitive element $x_0\in V_0$.
By Poincar\'e-Lefschetz duality
there is $x_1\in V_1$ such that $x_1\cap x_0=1$.
Since $H_3(C_1)=0$, there is $y\in H_4(W)$ such that $j_1y=x_1$.
We have $x_0\cdot j_0y=x_0\cap y=x_0\cap x_1=1$.
\qed

\smallskip
{\bf Bordism Theorem \Bd.}
{\it Let $(W,z)$ be a null-bordism of an admissible set
$$X=(C_0,C_1,A_0,A_1,\varphi)\quad\text{such that}\quad
\pi_1(C_k)=H_3(C_k)=H_4(C_k,\partial)=0. $$
The pair $(W,z)$ is spin bordant (relative to the boundary) to a product
with the interval if
there is a left inverse $s$ of the map
$$j:V_0\to H_4(W,\partial)$$
from the exact sequence of the triple $(W, \partial W, C_0)$, ($sj=\id$), such that }
$$\sigma(W)=sp_W\cdot sp_W=sz^2\cdot sp_W=sz^2\cdot sz^2=0.$$

{\it Beginning of the proof of the Bordism Theorem \Bd.}
Recall that $B\text{Spin}=BO\left<4\right>$ is the (unique up to homotopy) 3-connected
space for which there exists a fibration $B\text{Spin}\to BO$ inducing an isomorphism
on $\pi_i$ for $i\ge4$.
Denote $B :=B\text{Spin}\times\C P^\infty$.
Define $p:B\to BO$ to be the composition of the projection to $B\text{Spin}$ and
the map $B\text{Spin}\to BO$ inducing an isomorphism on $\pi_i$ for $i\ge4$.
Take the map $\overline\nu:W\to B$ corresponding to the given spin
structure on $W$ and to $z\in H_6(W,\partial)\cong[W,\C P^\infty]$.

Since $X$ is admissible and $H_4(C_k,\partial)=0$, by Poincar\'e-Lefschetz
duality the map $(\overline\nu|_{C_k})_*:H_2(C_k)\to H_2(\C P^\infty)$ is an
isomorphism.
This and $\pi_1(C_k)=0$ imply that the map $\overline\nu|_{C_k}$ is 3-connected.
Making $B$-surgery below the middle dimension we can change
$\overline\nu$ relative to the boundary and assume that $\overline\nu$ is
4-connected [Kr99, Proposition 4].
This surgery together with the obvious corresponding change of $s$ preserves
$\sigma(W),sp_W\cdot sp_W,sz^2\cdot sp_W$ and $sz^2\cdot sz^2$.
Hence it suffices to construct $U$ as in the [Kr99 Theorem \Kr].

Since $B\text{Spin}$ is 3-connected, we have
$$H_4(B)\cong H_4(B\text{Spin})\oplus H_4(\C P^\infty)\cong\Z\oplus\Z.$$
This isomorphism carries $\overline\nu_*u$ to $(u\cap p_W,u\cap z^2)$.
So `$\overline\nu_*U=0\in H_4(B)$' is equivalent to `$U\cap z^2=U\cap p_W=0$'.

Let
$$\widehat U=\{u\in V_0\ |\ au=msz^2+nsp_W\text{ for some integers }a,m,n\}.$$
(Note that $\rk \widehat U$ is 1 or 2.)
Since
$$sp_W\cdot sp_W=sz^2\cdot sp_W=sz^2\cdot sz^2=0,\quad\text{we have}
\quad \widehat U\cdot \widehat U=0.$$
Since the form $\cdot$ is unimodular, there is
$$X\subset V_0\quad \text{such that}\quad \widehat U\subset X,
\quad \rk X=2\rk \widehat U \quad \text{and \quad $\cdot|_X$ is unimodular}.$$
Then
\footnote{Since both $V_0$ and $X\subset V_0$ are unimodular, we have
$X\cap X^\perp=0$ and $\rk X^\perp=\rk V_0-\rk X$. Then $V_0=X\oplus X^\perp$.}
$V_0\cong X\oplus X^\perp$ and $\sigma(X)=0$.

The map $j_0:H_4(W)\to V_0$ is onto and carries $\cap$ to $\cdot$.
Therefore $\sigma(X^\perp)=\sigma(\cdot)=\sigma(W)=0$.
Hence there is a direct summand $\widetilde U\subset X^\perp$ such that
$\widetilde U\cdot \widetilde U=0$.
Let $U:=s^*(\widehat U\oplus \widetilde U)$, where $s^*$ is given by the
following Lemma \Se.

\smallskip
{\bf Lemma \Se.}
{\it Under the assumptions of Lemma \ho for each left inverse $s$ of $j$ a
right inverse $s^*:V_0\to H_4(W)$ of $j_0$ is well-defined by
$$s^*x\cap y=x\cdot sy\quad\text{for each}\quad y\in H_4(W,\partial).$$
The map $j_1s^*:V_0\to V_1$ is an
isomorphism carrying the product $\cap:V_0\times V_1\to\Z$ to $\cdot$,
i.e. $x\cdot x'=j_1s^*x\cap x'$ for each $x,x'\in V_0$.
\footnote{The second statement holds for each right inverse of $j_0$, not
necessarily the one obtained from $s$.}
}

\smallskip
{\it Proof.} Define a homomorphism $\overline x:H_4(W,\partial)\to\Z$
by $\overline x(y):=x\cdot sy$.
Now the existence and uniqueness of such an element $s^*x$ follows by
Poincar\'e-Lefschetz duality.

Clearly, $s^*$ is a homomorphism.

We have
$$j_0s^*x\cdot x'=s^*x\cap x'=s^*x\cap jx'=x\cdot sjx'=x\cdot x'
\quad\text{for each}\quad x,x'\in V_0.$$
Since the form $\cdot$ is unimodular, $j_0s^*x=x$.

We have $x\cdot x'=s^*x\cap x'=j_1s^*x\cap x'$.
(Cf. the end of the proof of Lemma \Ho.)


The map $s^*$ is injective.
For $x,x'\in V_0$ if
$$j_1s^*x=j_1s^*x',\quad\text{then}\quad x\cap a=j_1s^*x\cap a=j_1s^*y\cap a=
y\cap a\quad\text{for each}\quad a\in V_1.$$
Hence by Poincar\'e-Lefschetz duality
$x=y$.
Thus $j_1s^*$ is injective.
So it is an isomorphism.
\qed

\smallskip
{\it Completion of the proof of the Bordism Theorem \Bd:
checking the properties of $U$.}
Clearly, $\widehat U$ is a direct summand in $X$.

Let $U':=\widehat U\oplus \widetilde U$.
Then
$$j_0U=U',\quad U'\cdot U'=U'\cdot sz^2=U'\cdot sp_W=0$$
and $U'$ is a direct summand in $V_0$.

By Lemma \Se
$$U\cap U=U\cap jj_0U=s^*U'\cap jU'=U'\cdot sjU'=U'\cdot U'=0,$$
$$
U\cap x=U'\cdot sx=0\quad\text{for}\quad x\in\{z^2,p_W\}$$
and $j_0|_U$ is an isomorphism onto the direct summand $U'\subset V_0$.

Since $U\subset\im s^*$, by Lemma \se $j_1|_U$ is monomorphic.

Since $U'\subset V_0$ is a direct summand, we have $V_0\cong U'\oplus U''$
(additive) for some $U''\subset V_0$.
Suppose that $j_1s^*u'=j_1s^*u''$ for some $u'\in U'$ and $u''\in U''$.
By excision $H_4(\partial W,C_1)\cong H_4(C_0,\partial)=0$, so by the exact
sequence of pair the inclusion-induced map $H_4(C_1)\to H_4(\partial W)$ is
surjective. Hence for the inclusion-induced maps
$$i:H_4(\partial W)\to H_4(W)\quad\text{and}\quad i_k:H_4(C_k)\to H_4(W)\quad
\text{we have}\quad\im i=\im i_1.$$
Analogously $\im i=\im i_0$.
Hence
$$s^*u'-s^*u''\in\im i_1=\im i_0,
\quad\text{so}\quad u'-u''=j_0(s^*u'-s^*u'')=0,
\quad\text{hence}\quad u'=u''=0.$$
Thus $j_1U\cap j_1s^*U''=0$.
Therefore by dimension considerations $V_1\cong j_1U\oplus j_1s^*U''$
(additively).  So $j_1U$ is a direct summand.

The pairing $\cap:j_0U\times V_1/j_1U\to\Z$ is isomorphic to the pairing
$\cap:U'\times j_1s^*U''\to\Z$ and (by Lemma \Se) to the pairing
$\cdot:U'\times U''\to\Z$.
Since the form $\cdot:V_0\times V_0\to\Z$ is unimodular and $U'\cdot U'=0$, the
latter pairing is unimodular.
\qed

\bigskip
{\bf Proof of the `if' part of the Almost Diffeomorphism Theorem \Hs.}

\smallskip
{\it Beginning of the proof.}
Take a null-bordism $(W,z)$ of $X$ given by the Null-bordism Lemma \Co.
The idea is to modify $(W,z)$ and apply the Bordism Theorem \Bd.
Define $B,p$ and a 4-connected map $\overline\nu:W\to B$ as in the beginning of
the proof of the Bordism Theorem \Bd.

Since $H_3(C_0)=0$, we can take the product $\cdot$ given by Lemma \Ho.

By excision $H_4(\partial W,C_0)\cong H_4(C_1,\partial)=0$.
Then, by the exact sequence of a triple, the homomorphism $j : V_0 \to H_4(W, \partial)$ is injective.

Take $x\in V_0$.
We have $x'\cdot x=y\cap x=y\cap jx$ for each $x'\in V_0$ and
$y\in j_0^{-1}x'$.
If $jx$ is divisible by an integer $d$, then $x'\cdot x$ is divisible by $d$
for each $x'\in V_0$.
Hence the unimodularity of $\cdot$ implies
that $jx$ is primitive for each primitive $x\in V_0$.
So there exists a left inverse $s$ of $j$ (because $\overline\nu$
is 4-connected and so $\Tors H_4(W,\partial)=\Tors H_3(W)=0$).

Denote $d:=d(\partial_W z^2)$.
Recall the definition of $\overline{p_W}\in H_4(W)$ and
$\overline{z^2}\in H_4(W;\Z_d)$ from the definition of $\eta_X$ in \S2.
Since $j_0\overline{p_W}=sp_W$, we have $\overline{p_W}\cap p_W=sp_W\cdot sp_W$.
Since
$$j_0\overline{z^2}=\rho_dsz^2,\quad\text{we have}\quad
\overline{z^2}\cap p_W=\rho_dsz^2\cdot sp_W\in\Z_d\quad\text{and}\quad
\overline{z^2}\cap z^2=\rho_dsz^2\cdot sz^2\in\Z_d.$$
Denote $\widehat\eta_{W, z, s}=sz^2\cdot(sz^2-sp_W)\in\Z$.
Thus $\eta_X=\rho_d\widehat\eta_{W, z, s}$.
\footnote{Note that $\rho_d(\overline{p_W}\cap z^2)=\overline{z^2}\cap\rho_dp_W
=\rho_d(sp_W\cdot sz^2)$ but
$\overline{p_W}\cap z^2\ne sp_W\cdot sz^2=
s^*j_0\overline{p_W}\cap z^2$.}

Analogously if $A_0^2$ is divisible by 2 then $\eta'_X=\rho_2(sz^2\cdot sz^2)$.

\smallskip
For the completion of the proof we need two lemmas.
Let $W$ be a compact spin 8-manifold such that $\partial_Wp_W=0$.
Then there is $\overline{p_W}\in H_4(W)$ such that $j_W\overline{p_W}=p_W$.
(It is clear that the intersections below do not depend on the choice of
$\overline{p_W}$.)
By Lemma \Di(b)
$$\sigma(W)
\underset{\mod8}\to\equiv \overline{p_W}\cap p_W\quad\text{so}
\quad \alpha_W:=\frac{\sigma(W)-\overline{p_W}\cap p_W}8\quad
\text{is an integer}.$$

{\bf Lemma \Stw.} {\it For each of the four quadruples
$$(1,0,0,0),\quad(0,28,0,0),\quad(0,0,2,0),\quad(0,0,0,12)$$
there is a closed compact spin 8-manifold $W$ and
$z\in H_6(W)$ such that the quadruple
$Q_{W,z}:=(\sigma(W),\alpha_W,z^4,
\ z^4-z^2p_W)$
coincides with the given quadruple.
\footnote{We could avoid using $(0,0,2,0)$ by using the Framing Theorem
\Ind($\varphi$) and changing the structure of the proof of the injectivity of
$\eta_u$.}
}

\smallskip
{\bf Lemma \Tw.} {\it
Let $(W,z)$ be a null-bordism of an admissible set $X$ such that $H_3(C_k)=H_3(W)=H_5(W,\partial)=0$.
Let $s$ be a left inverse of $j$.
By connected sum of $W$ with a null-bordant closed 3-connected 8-manifold and
certain corresponding change of $z,s$ one can change:

$\bullet$ $sz^2\cdot sz^2$ by adding an odd number, provided $A_0^2$ is not divisible
by 2,

$\bullet$ $\widehat\eta_{W, z, s}$ by adding $2d/GCD(d,2)$, where $d:=d(A_0^2)$,
and preserving $\rho_2(sz^2\cdot sz^2)$.}

\smallskip
The lemmas are proved in the next subsection (Lemma \stw is known).

\smallskip
{\it Completion of the proof of the `if' part of the Almost Diffeomorphism
Theorem \Hs.}
Take a 3-connected parallelizable $8$-manifold $\overline E_8$ whose boundary
is a homotopy sphere and whose signature is $8$.
Then $p_{\overline E_8}=0$.
The boundary connected sum of $\overline\nu$ with a constant map
$\overline E_8\to \C P^\infty$ changes $\alpha_W$
by 1 and preserves the 4-connectedness of $\overline\nu$.
\footnote{An alternative proof is obtained by replacing $\overline E_8$ by
a $3$-connected $8$-manifold $X\simeq S^4$ whose boundary is
a homotopy sphere, $\sigma(X)=1$ and $p_X=3$ [Mi56].}
Thus we may assume that $\alpha_W=0$.


For a null-bordism $(W,z)$ of an admissible set $X$ such that $H_3(C_k)=0$
and a left inverse $s$ of $j$
denote $$Q_{(W, z, s)}:=(\sigma(W),\alpha_W,sz^2\cdot sz^2,\widehat\eta_{(W, z, s)}).$$
For a closed spin 8-manifold $W_0$ and $z_0\in H_6(W_0)$
we have $$Q_{W\#W_0,z\oplus z_0,s\oplus\id}=Q_{(W, z, s)}+Q_{W_0,z_0}.$$
Since $z$ is primitive, $z\oplus z_0$ is primitive.
So we may spin surger $W\#W_0$ and assume that the map
$\overline\nu':W\#W_0\to B$ corresponding to $z\oplus z_0$ and the
`connected sum' spin structure on $W\#W_0$ is 4-connected.
So by Lemma \stw we may change the quadruple $Q_{W, z, s}$ by any of the four
quadruples of Lemma \Stw, and $\overline\nu$ would remain 4-connected.

Thus we may assume that $\sigma(W)=\alpha_W=0$.

Taking the connected sum of $\overline\nu$ with the constant map from a null-bordant
3-connected 8-manifold does not change $\sigma(W)$, $\alpha_W$, or the
property that $\overline\nu$ is 4-connected.

If $A_0^2$ is not divisible by 2, then by Lemmas \tw
and \stw we may assume that $\sigma(W)=\alpha_W=sz^2\cdot sz^2=0$.

If $A_0^2$ is divisible by 2, then $\rho_2(sz^2\cdot sz^2)=\eta'_X=0$, hence
by Lemma \stw we may assume that
$\sigma(W)=\alpha_W=sz^2\cdot sz^2=0$.

Since $\eta_X=0$, by Lemmas \tw and \stw we may assume
that $\sigma(W)=\alpha_W=sz^2\cdot sz^2=\widehat\eta_{(W, z, s)}=0$.
Then we are done by the Bordism Theorem \Bd.
\qed

\smallskip
{\bf Diffeomorphism Theorem \Dif.}
{\it Let $X=(C_0,C_1,A_0,A_1,\varphi)$ be an admissible set such that
$\pi_1(C_k)=H_3(C_k)=H_4(C_k,\partial)=0$.
Denote $\alpha_X:=\rho_{28}\alpha_W\in\Z_{28}$ for some null-bordism $(W,z)$ of $X$.
\footnote{The independence of $\alpha_X$ from $W$ is easily deduced from known results.
Note that $\alpha_X$ is also independent of $\varphi$ because for a closed
spin 8-manifold $V$ we have
$\sigma(V)-p_V^2\in (2^5 \cdot 7) \cdot \Z$.
These remarks are not necessary for our main results.
}
There is a diffeomorphism $\overline\varphi:C_0\to C_1$ extending $\varphi$ and such that $\overline\varphi_*A_0=A_1$ if and only if
$$\alpha_X=0,\quad\eta_X=0\quad\text{and, for $A_0^2$ divisible by 2,}\quad
\eta'_X=0.$$}
\quad
The `only if' part is simple and is proved analogously to the `only if' part of the Almost Diffeomorphism Theorem \hs (we do not need $V$ and define $W:=C_0\times I$, $z:=A_0\times I$).
The proof of the `if' part is very close to the proof of the `if' part of the Almost Diffeomorphism Theorem \Hs.
The only change is that in the completion of the proof of the `if' part we have $\alpha_W=0$
by hypothesis and do not use the boundary connected sum of $\overline\nu$ with the constant map $\overline E_8\to\C P^\infty$

\smallskip
{\bf Bordism Conjecture \Cjc.}
{\it Let $(W,z)$ be a null-bordism of an admissible set
$X=(C_0,C_1,A_0,A_1,\varphi)$ such that
$\pi_1(C_k)=H_3(C_k)=H_4(C_k,\partial)=0$,
and the map $h_z:W\to\C P^\infty$ corresponding to $z$ is 4-connected.
Then $(W,z)$ is spin bordant (relative to the boundary) to a product
with the interval if and only if}
$$\sigma(W)=\overline{p_W}\cap p_W=0\quad\text{and}\quad
\overline{z^2}\cap p_W=\overline{z^2}\cap z^2=0\in\Z_d.$$

\smallskip
{\bf Proof of Lemmas \stw and \Tw.}

Denote $H_+:=\left(\matrix0&1\\1&0\endmatrix\right)$.

\smallskip
{\it Proof of Lemma \Stw.}
Recall that $\sigma(\Ha P^2)=1$ and $p_1^2(\Ha P^2)=4$ [Hi53], cf.
[Mi56, Lemmas 3 and 4].
Thus by Lemma \Di.a for $(\Ha P^2,0)$ the quadruple is $(1,0,0,0)$.

Take a 3-connected parallelizable 8-manifold $\overline E_8$ whose boundary is
a homotopy sphere and whose signature is 8.
Then $p_1(\overline E_8)=0$.
Thus by Lemma \Di.a for $(28\overline E_8\cup D^8,0)$ the quadruple is
$(28\cdot8,28,0,0)$.

Take $(S^2)^4$ and the class $z$ which is the sum of four summands, each
represented by a product of three 2-spheres and a point.
Then $z^4=24$.
As a quadratic form $H_4((S^2)^4)\cong H_+\oplus H_+\oplus H_+$, so
$\sigma((S^2)^4)=0$.
Since $(S^2)^4$ is almost parallelizable, we have $p_1((S^2)^4)=0$.
Thus by Lemma \Di.a for $((S^2)^4,z)$ the quadruple is $(0,0,24,24)$.

By [KS91, Proposition 2.5] there is a closed spin 8-manifold $W$ and
$z\in H_6(W)$ such that $S_1=S_2=0$ and $S_3=1$.
In the notation of [KS91, spin case of (2.4)]
$$S_1=\alpha_W/28,\quad S_2=z^2(z^2-p_W)/12\quad\text{and}
\quad 2S_3=8S_2+z^4.$$
Hence for $(W,z)$ the quadruple is $(a,0,2,0)$.
\qed

\smallskip
{\bf Lemma \Sh.}
{\it Assume that $(W,z)$ is a null-bordism of an admissible set $X$.}

(p) {\it $s'p_W=sp_W$ for each left inverses $s,s'$ of $j$.}

(z) {\it Suppose that $H_3(C_0)=H_3(W)=H_5(W,\partial)=0$.
For $x\in V_0$ there is a left inverse $s'$ of $j$ such that
$s'z^2=sz^2+x$ if and only if $x$ is divisible by $d:=d(\partial_W z^2)$.}

\smallskip
{\it Proof of }(p).
Denote by $\partial_0:H_4(W,\partial)\to H_3(\partial W,C_0)$ the boundary
homomorphism.
The class $(\partial p_W)\cap C_0=PDp_1(C_0)=0$ goes to $\partial_0 p_W$ under
the excision isomorphism $H_3(C_1,\partial)\to H_3(\partial W,C_0)$.
Thus $\partial_0 p_W=0$.
Hence $p_W\in\im j$ which implies (p).
\qed

\smallskip
{\it Proof of }(z).
Since $H_3(C_0)=0$, the map $j_0$ is onto, hence
$\im j=\im(jj_0)=\ker\partial_W$.
Since $H_3(\partial C_0)=0$, we have that $H_2(\partial C_0)$ is torsion free. Using this and  $H_3(C_k)=0$, by the Mayer-Vietoris
sequence for $\partial W=C_0\cup C_1$ we obtain that $H_3(\partial W)$ is torsion free.
This and $H_3(W)=H_5(W,\partial)=0$ imply that
$H_4(W,\partial)\cong V_0\oplus H_3(\partial W)$.
Identify these isomorphic groups by the isomorphism
$j\oplus(\partial_W|_{\ker s})^{-1}$.
Then $z^2$ is identified with $sz^2\oplus\partial_Wz^2$.
The `only if' part follows because $s'(sz^2\oplus0)=sz^2$, so
$s'z^2=sz^2+s'\partial_Wz^2$.
The `if' part follows because
$\partial_Wz^2/d\in H_3(\partial W)\subset H_4(W,\partial)$ is
primitive, so for each $x_1\in V_0$ there is a left inverse $s'$ of $j$ such
that $s'(z^2/d)=s(z^2/d)+x_1$.
\qed

\smallskip
{\it Proof of Lemma \Tw.}
First we prove the second assertion.
By [Mi56] there is a $D^4$-bundle over $S^4$ whose Euler class is 0
and whose first Pontryagin class is 4.
The double of this bundle is an $S^4$-bundle $S^4\widetilde\times S^4$ over
$S^4$ whose first Pontryagin class is $4$.
We have $H_4(S^4\widetilde\times S^4)\cong\Z\oplus\Z$ with evident basis.
In this basis $p_1(S^4\widetilde\times S^4)=(4,0)$ and the intersection form of
$S^4\widetilde\times S^4$ is $H_+$.

Denote $W':=W\#(S^4\widetilde\times S^4)$.
Identify $H_6(W,\partial)$ with $H_6(W',\partial)$.
Identify $H_4(W',C_0)$ with $V_0\oplus H_+$ as groups with quadratic forms.
Clearly,
$$\partial W'=\partial W,\quad\partial_{W'}z=\partial_{W}z
\quad\text{and}\quad \widehat\eta_{W',z,s\oplus\id}=\widehat\eta_{W,z,s}.$$
By Lemma \bo and (the `if' part of) Lemma \Sh(z) there is a left inverse
$$s':H_4(W',\partial)\to  H_4(W',C_0)\quad\text{such that}
\quad s'(z^2\oplus(0,0))=sz^2\oplus(0,d).$$
By Lemma \Di(a) $p_{W'}=p_W\oplus(2,0)$.
By Lemma \Sh(p), $s'p_{W'}=(s\oplus\id)p_{W'}=sp_W\oplus(2,0)$.
So
$$sz^2\cdot sz^2=s'z^2\cdot s'z^2\quad\text{and}\quad
\eta_{W',z,s'}-\eta_{W,z,s}=(0,d)\cap[(0,d)-(2,0)]=(0,d)\cap(-2,d)=-2d.$$
\quad
In this paragraph assume that $d$ is even.
We have $H_4(\Ha P^2\#(-\Ha P^2))\cong\Z\oplus\Z$ with evident basis.
In this basis $p_1(\Ha P^2\#(-\Ha P^2))=(2,-2)$ and the intersection form of
$\Ha P^2\#(-\Ha P^2)$ is $\diag(1,-1)$.
Analogously to the above with $S^4\widetilde\times S^4$ replaced by
$\Ha P^2\#(-\Ha P^2)$ we may change $\eta_{(W, z, s)}$ by
$$(0,d)\cap[(0,d)-(1,-1)]=(0,d)\cap(-1,d+1)=-d^2-d.$$
The difference $s'z^2\cdot s'z^2-sz^2\cdot sz^2=(0,d)\cap(0,d)=-d^2$ is
divisible by 2.
Hence we may change $\widehat\eta_{W,z,s}$ by $GCD(2d,d^2+d)=d$ and preserve
$\rho_2(sz^2\cdot sz^2)$.

Now we prove the first assertion.
Since $A_0^2$ is not divisible by 2, $d$ is odd.
Hence in the above argument involving $\Ha P^2\#(-\Ha P^2)$ the change of
$sz^2\cdot sz^2$ is by an odd integer $d^2$.
\qed

\head 5. Remarks  \endhead

It would be interesting to know which embeddings $f:N\to\R^7$ of
a closed orientable 4-manifold $N$ have Seifert surfaces.

The following properties from the definition of admissibility are
not necessary for some lemmas:

$H_5(\partial C_0)=0$ is only used for the uniqueness of $\partial z$;

$H_3(\partial C_0)=0,\quad p_{C_0}=p_{C_1}=0\quad\text{and}
\quad d(A_0^2)=d(A_1^2)$ for the Null-bordism Lemma \Co,

$d(A_0^2)=d(A_1^2)$ for the definition of $\eta_{(W, z)}$,

$d(A_0^2)=d(A_1^2)$ and $p_{C_0}=p_{C_1}=0$ for the definition of $\eta'_X$
and the Bordism Theorem \Bd,

$p_{C_0}=p_{C_1}=0$ for the Framing Theorem \Ind,

$d(A_0^2)=d(A_1^2)$ and $H_3(\partial C_0)=0$ for Lemmas \tw  and \Sh.

\smallskip
{\it Remarks to the construction of a 1--1 correspondence between normal framings
on an embedding $S^3\to S^7$ (up to homotopy) and $\Z\oplus\Z$.}
Surgery on a framed embedding $b:S^3\times D^4\to S^7$ gives a 8-manifold
$E_b$ which is the total space of a $D^4$-bundle $E_b\to S^4$.
The boundary $\partial E_b$ is the total space of an $S^3$-bundle $\xi_b:E_b\to S^4$. The map $b\mapsto\xi_b$ is a 1--1 correspondence [Wa62, Lemma 1].
Take the 1--1 correspondence between $S^3$-bundles over $S^4$ and $\Z\oplus\Z$
constructed in [Mi56].
This gives an alternative construction of the above 1--1 correspondence.

The map assigning to $b$ the diffeomorphism class of the total space
$E_b$ is a bijection.
The inverse is given by
$E\mapsto(\dfrac{a_E\cap a_E-p_E\cap a_E}2,\dfrac{a_E\cap a_E+p_E\cap a_E}2)$,
where $a_E\in H_4(E)$ is the generator and we use the above 1--1 correspondence
between the set of framings and $\Z\oplus\Z$.
\footnote{
The map assigning to $b$ the diffeomorphism class of the total space
$\partial E_b$ is {\it not} a bijection (although the restriction of such a map
gives a 1--1 correspondence between {\it unlinked} framed embeddings and
diffeomorphism classes of total spaces of trivial Euler class bundles) [CE03].
\newline
Framed embeddings $b$ corresponding to pairs $(a,-a)$ are
characterized by being {\it unlinked} (i.e. such that the linking coefficient
of $b(S^3\times 0)$ and $b(S^3\times x)$ is zero.
\newline
An isotopy $F$ from an embedding $S^3\to S^7$ to the standard
embedding is not necessarily unique up to isotopy (of isotopies relative to the
ends).
So apriori we cannot just take as the 'zero' framing the image of
the standard framing of the standard embedding under such an isotopy $F$.
However, the above argument shows that we can.
}

\bigskip
{\bf An alternative proof of the Agreement Lemma.}

The Agreement Lemma is an analogue of [Sk08', the Agreement Lemma].
For $H_1(N)\ne0$ this analogue is more complicated
 because embeddings $N_0\to S^7$ are not necessarily isotopic.

Let $f:N\to S^7$ be an embedding.
In this subsection we omit subscript $f$ of $\nu_f,C_f,A_f$ etc.
A section $\xi:N_0\to\nu^{-1}N_0$ is called {\it faithful} if
$\xi^!\partial A=0$.
When $H_2(N)$ has no torsion, this is equivalent to the triviality of the
composition
$H_2(N_0)\overset{\xi_*}\to\to H_2(\partial C)\overset{i_*}\to
\to H_2(C)$.

Faithfulness is not equivalent to unlinkedness because in general
$AD_{f|_{N_0}}\overline\xi_* \ne f|_{N_0}^!AD_{\overline\xi}$.

The Agreement Lemma is implied by the following result.

\smallskip
{\bf Faithful Section Lemma.}
(a) {\it A faithful section exists.
It is unique on 2-skeleton of $N$ up to fiberwise homotopy.
[HH63, 4.3, BH70, Proposition 1.3].}

(b) {\it Under the assumptions of the Agreement Lemma $\varphi$ maps
a faithful section to a faithful section.}

\smallskip
Part (a) is implied by the following result.

\smallskip
{\bf Difference Lemma.}
{\it For each pair of sections $\xi,\eta:N_0\to\nu^{-1}N_0$ we have
\linebreak
$d(\xi,\eta)=(\xi^!-\eta^!)ej\partial A$.}

\smallskip
Let $\zeta$ be an unlinked section for $f$.

\smallskip
{\it Sketch of a proof of the Difference Lemma.}
The lemma follows because
$(\xi^!-\eta^!)ej\partial A=(\xi^!-\eta^!)\zeta_*[N_0]=d(\xi,\eta)$.
Here the first equality holds by [Sk10, Section Lemma 2.5.a]
because $\zeta$ is unlinked.
Let us sketch a proof the second equality
(for any section $\zeta$; this equality generalizes to any bundle;
note that in general $\xi^!\zeta_*[N_0]\ne d(\xi,\zeta)$.)
Recall that $d(\xi,\eta)$ is the intersection in $\widehat\nu^{-1}N_0$
of $fN_0$ and a 5-film $\Delta$ with spanned by $\xi N_0$ and $\eta N_0$.
Take a 5-film in $\widehat\nu^{-1}N_0$ spanned by $fN_0$ and $\zeta N_0$.
We may assume that these two 5-films are in general position to each other
so that their intersection is a homology between $d(\xi,\eta)$ and
$[\Delta\cap\zeta N_0]=(\xi^!-\eta^!)\zeta_*[N_0]$.
\qed

\smallskip
{\it Proof of the Faithful Section Lemma} (b).
Recall the equality on $\pm2d(\xi,\eta)$
from the proof of the Agreement Lemma in \S3.
Then for a faithful section $\xi$ for $f$ we have
$$PDe(\zeta^\perp)-PDe(\xi^\perp)=\pm2d(\zeta,\xi)=
\pm2(\zeta^!-\xi^!)\partial A_f=\pm2\zeta^!\partial A_f=\pm2PDe(\zeta^\perp).$$
Here

$\bullet$ the first equality holds by  the equality on $\pm2d(\xi,\eta)$;

$\bullet$ the second equality holds by the Difference Lemma,

$\bullet$ the third equality holds because $\xi$ is faithful,

$\bullet$ the fourth equality holds by (the second equality of) the Section Lemma.

Since $H_2(N)$ has no 2-torsion, together with the equality
on $\pm2d(\xi,\eta)$  this implies that either

{\it a section $\xi:N_0\to\partial C_f$ is faithful if
and only if $PDe(\xi^\perp)=-PDe(\zeta^\perp)$}, or

{\it a section $\xi:N_0\to\partial C_f$ is faithful if
and only if $PDe(\xi^\perp)=3PDe(\zeta^\perp)$.}

Now the lemma follows by the Section Lemma \eu because
$e((\varphi\xi)^\perp)=e(\xi^\perp)$.
\qed

\smallskip
We conjecture that $\varkappa (f)-\varkappa (f')=2W_{f'}(f)$ for the {\it Whitney
invariant} $W_{f'}(f)$ [Sk08, \S2].
For simply-connected $N$ the proof is analogous to [Sk08', \S3].

The following assertion is proved analogously to [Sk08', the Difference Lemma
(c)] (where $A_0$ is defined):
{\it if $f=f'$ on $N_0$ and $\xi:N_0\to\partial C_f$ is a section both for
$f$ and $f'$, then $W(f)-W(f')=A_0(\overline\xi_*-\overline\xi'_*)[N]$, where
$\overline\xi'$ is constructed from $\xi$ and $f'$.}

This assertion gives an alternative proof
of the following statement used in the proof of the Agreement Lemma:
{\it if $\varkappa (f)=\varkappa (f')$ and $H_1(N)=0$, then any bundle isomorphism maps an unlinked
section of $f$ to that of $f'$.}
\footnote{If a section $\xi:N_0\to\partial C_f$ is strongly unlinked, then
it is faithful.
If $N$ is simply-connected, then the converse also holds because
$N_0\simeq\vee S^2_i$.
If a section $\xi:N\to\partial C_f$ is strongly unlinked, then its restriction
to $N_0$ is both faithful and unlinked, hence $\varkappa (f)=0$ by the italicized
assertion in the proof of the Faithful Section Lemma (b).
The same assertion implies that
{\it for simply-connected $N$ the existence of a strongly unlinked
framing of $\nu_0$ is equivalent to $\varkappa (f)=0$ (and hence to the
compressibility of $f$).}
Here the simply-connectedness assumption is essential:
take an embedding $(S^1\times S^3)_1\#(S^1\times S^3)_2$
such that $(x\times S^3)_1$ and $(x\times S^3)_2$ are linked, then for any
section $\xi:N_0\to\partial C_f$ we have
$\xi^*in^*\ne0\in H^3(N_0)$.
If $\nu$ is trivial, then the obstruction to extending a section
$\xi:N_0\to\partial C_f$ to $N$ is $(\xi^!\partial A_f)^2\in\Z$.
Thus unlinked or faithful section on $N_0$ extends to $N$ if and only if
$\varkappa (f)=0$.)
}


\enddocument